\documentclass[reqno,11pt]{amsart}
\usepackage{amssymb}
\usepackage{amsmath}
\usepackage{latexsym,array}
\usepackage{amsfonts}
\newcommand{\C}{\mathbb {C}}

\newcommand{\bp}{\boldsymbol{\rho}}
\newcommand{\balpha}{\boldsymbol{\alpha}}

\newcommand{\sbm}[1]{\left[\begin{smallmatrix} #1
                \end{smallmatrix}\right]}

\newcommand{\cX}{\mathcal {X}}

\newcommand{\bH}{\mathbb{H}}

\newcommand{\R}{\mathbb{R}}

\newcommand{\bl}{\boldsymbol{e_\ell}}
\newcommand{\br}{\boldsymbol{e_r}}

\newtheorem{Pa}{Paper}[section]
\newtheorem{theorem}[Pa]{{\bf Theorem}}
\newtheorem{lemma}[Pa]{{\bf Lemma}}
\newtheorem{definition}[Pa]{{\bf Definition}}
\newtheorem{corollary}[Pa]{{\bf Corollary}}

\newtheorem{remark}[Pa]{{\bf Remark}}
\newtheorem{proposition}[Pa]{{\bf Proposition}}

\newtheorem{example}[Pa]{{\bf Example}}

\begin{document}
\title[]{Confluent Vandermonde matrices, divided differences, and Lagrange-Hermite
interpolation over quaternions}
\author{Vladimir Bolotnikov}
\email{vladi@math.wm.edu}
\address{Department of Mathematics,
The College of William and Mary,
Williamsburg VA 23187-8795, USA}
\begin{abstract}
We introduce the notion of a confluent Vandermonde matrix with quaternion entries
and discuss its connection with Lagrange-Hermite interpolation over quaternions.
Further results include the formula for the rank of a confluent Vandermonde matrix,
the representation formula for divided differences of quaternion polynomials and their 
extensions to the formal power series setting.
\end{abstract}
\maketitle

\section{Introduction}
\setcounter{equation}{0}

The notion of the Vandermonde matrix arises naturally in the context of 
the Lagrange interpolation problem when one seeks a complex 
polynomial taking prescribed values at given points. Confluent Vandermonde matrices 
come up once interpolation conditions on the derivatives of an unknown interpolant are also 
imposed; we refer to the survey \cite{kalman} for confluent Vandermonde matrices
and their applications. The study of Vandermonde matrices over division rings (with further 
specialization to the ring of quaternions) was initiated in \cite{lam1}. In the follow-up 
paper \cite{lamler}, the Vandermonde matrices were studied in the setting of a division ring 
$K$ endowed with an endomorphism $s$ and an $s$-derivation $D$, along with their interactions
with skew polynomials from the Ore domain $K[z,s,D]$. The objective of this paper is to extend 
the results from \cite{lam1} in a different direction: to introduce a meaningful
notion of a {\em confluent Vandermonde matrix} (over quaternions only, for the sake of 
simplicity) and to discuss its connections with Lagrange-Hermite interpolation problem for 
quaternion polynomials. 

\smallskip

Let $\bH$ denote the skew field of quaternions
$\alpha=x_0+{\bf i}x_1+{\bf j}x_2+{\bf k}x_3$ where  $x_0,x_1,x_2,x_3\in\mathbb R$   
and where ${\bf i}, {\bf j}, {\bf k}$ are the imaginary units commuting with $\R$ and
satisfying ${\bf i}^2={\bf j}^2={\bf k}^2={\bf ijk}=-1$. For $\alpha\in\bH$ as above,
its real and imaginary parts, the quaternion conjugate and the absolute value
are defined as ${\rm Re}(\alpha)=x_0$, ${\rm Im}(\alpha)={\bf i}x_1+{\bf j}x_2+{\bf k}x_3$,
$\overline \alpha={\rm Re}(\alpha)-{\rm Im}(\alpha)$ and
$|\alpha|^2=\alpha\overline{\alpha}=|{\rm Re}(\alpha)|^2+|{\rm Im}(\alpha)|^2$,
respectively. Two quaternions $\alpha$ and $\beta$ are called {\em equivalent} 
(conjugate to each other) if $\alpha=h^{-1}\beta h$ for some nonzero $h\in\mathbb H$; in 
notation, $\alpha\sim\beta$. It turns out that
\begin{equation}
\alpha\sim\beta\quad\mbox{if and only if}\quad {\rm Re}(\alpha) ={\rm Re}(\beta) \;
\mbox{and} \; |\alpha|=|\beta|,
\label{1.0}   
\end{equation}
so that the {\em conjugacy class} of a given $\alpha\in\mathbb H$ form a $2$-sphere (of
radius $|{\rm Im}(\alpha)|$ around ${\rm Re}(\alpha)$) which will be denoted  by $[\alpha]$.
It is clear that $[\alpha]=\{\alpha\}$ if and only if $\alpha\in\R$. For an 
$\alpha\in\bH\backslash\R$,
we denote by $\C_\alpha$ the plane spanned by $\alpha$ and $1$, which alternatively can be 
characterized as the set of all quaternions commuting with $\alpha$. Observe that $\C_\alpha\cap[\alpha]=
\{\alpha,\overline{\alpha}\}$.

\medskip

\begin{definition}
{\rm A finite ordered collection ${\balpha}=(\alpha_1,\ldots,\alpha_k)$ is a {\em 
spherical chain} (of length $k$) if
\begin{equation}
\alpha_1\sim \alpha_2\sim\ldots\sim\alpha_k\quad\mbox{and}\quad \alpha_{j+1}\neq
\overline{\alpha}_j\quad\mbox{for}\quad j=1,\ldots,k-1.
\label{1.8}
\end{equation}}
\label{D:1.3}   
\end{definition}
The latter notion is essentially non-commutative: a spherical chain $\balpha$ consisting of
commuting elements is necessarily a subset of $\C_{\alpha_1}\cap [\alpha_1]\subset \{\alpha_1,\overline{\alpha}_1\}$ which
together with inequality in \eqref{1.8} implies that all elements in $\balpha$ are the same:
\begin{equation}
\balpha=(\alpha,\alpha,\ldots,\alpha),\qquad \alpha\in\bH.
\label{1.9}
\end{equation}
The spherical chain \eqref{1.9} may be understood as an element $\alpha$ taken with
multiplicity $k$. Observe that many (if not all) results concerning multiple zeros of complex polynomials
or multiple eigenvalues of complex matrices extend to the quaternion setting almost literally
and can be interpreted in terms of commutative chains of the form \eqref{1.9}. It turns out, however, 
that 
they may have meaningful and nontrivial extensions to more general (non-commutative) spherical chains 
\eqref{1.8}.
The latter observation is partly justified by the existing results on the zero structure of quaternion polynomials
\cite{bol2,genstr1,gm} and canonical forms for quaternion matrices \cite{brenner,lee,wieg}; some 
extra 
evidence (the structure of the confluent Vandermonde matrices and interpolation conditions of Lagrange-Hermite type 
which do not appear in the commutative case) will be given below. 

\bigskip
\noindent
{\bf 1.1. The basic case.} We now briefly review the basic (non-confluent) results  in the form convenient for further 
``confluent" extensions. Let $\bH[z]$ denote the ring of polynomials in one formal variable $z$ which commutes with
quaternion coefficients. The ring operations in $\bH[z]$ are defined as in the commutative
case, but as multiplication in $\bH$ is not commutative, multiplication in $\bH[z]$ is not 
commutative either. A straightforward computation verifies that for any $\alpha\in\bH$ and 
$f\in\bH[z]$,
\begin{equation}
f(z)=f^{\bl}(\alpha)+(z-\alpha)\cdot(L_\alpha f)(z)=f^{\br}(\alpha)+(R_\alpha f)(z)\cdot(z-\alpha),
\label{1.1}
\end{equation}
where $f^{\bl}(\alpha)$ and $f^{\br}(\alpha)$ are respectively, left and right evaluation of
$f$ at $\alpha$ given by
\begin{equation}
f^{\bl}(\alpha)=\sum_{k=0}^m\alpha^k f_k\quad\mbox{and}\quad
f^{\br}(\alpha)=\sum_{k=0}^m f_k\alpha^k\quad\mbox{if}\quad f(z)=\sum_{j=0}^m z^j f_j,
\label{1.2}
\end{equation}
and where $L_\alpha f$ and $R_\alpha f$ are polynomials of degree $m-1$ given by
\begin{equation}(L_\alpha f)(z)={\displaystyle\sum_{k=0}^{m-1}\bigg(\sum_{j=0}^{m-k-1}
\alpha^jf_{k+j+1}\bigg)z^k},\qquad
(R_\alpha f)(z)={\displaystyle\sum_{k=0}^{m-1}\bigg(\sum_{j=0}^{m-k-1}
f_{k+j+1}\alpha^j\bigg)z^k}.
\label{1.3}
\end{equation}
Interpreting $\bH[z]$ as a vector space over $\bH$, we observe that  
the mappings $f\mapsto L_\alpha f$ and $f\mapsto R_\alpha f$
define respectively the right linear operator $L_\alpha$ and the left
linear operator $R_\alpha$ (called in analogy to the complex case, the left and the 
right backward shift, respectively) acting on $\bH[z]$.

\smallskip

With evaluations \eqref{1.2} in hands, we formulate the left Lagrange 
interpolation problem in $\bH[z]$: {\em given distinct ``points" 
$\alpha_1,\ldots,\alpha_n\in\bH$ and
target values $c_1,\ldots,c_n\in\bH$, find a polynomial $f\in\bH[z]$
such that
\begin{equation}
f^{\bl}(\alpha_i)=c_i\quad\mbox{for}\quad i=1,\ldots,n,
\label{1.4}
\end{equation}}
or the right interpolation problem: {\em find a polynomial $f\in\bH[z]$
such that
\begin{equation}
f^{\br}(\alpha_i)=c_i\quad\mbox{for}\quad i=1,\ldots,n.
\label{1.5}   
\end{equation}}
Making use of \eqref{1.2} (with $m-1$ instead of $m$) we can write conditions 
\eqref{1.4}, \eqref{1.5} in the form
\begin{align*}
f_0+\alpha_if_1+\ldots \alpha_i^{m-1}f_{m-1}&=c_i\qquad (i=1,\ldots,n),\\
f_0+f_1\alpha_i+\ldots f_{m-1}\alpha_i^{m-1}&=c_i\qquad (i=1,\ldots,n),
\end{align*}
respectively, with unknown $f_0,\ldots,f_{m-1}$, or in the matrix form, as
$$
V_m^{\boldsymbol\ell}F=C,\quad F^{\top}V_m^{\boldsymbol r}=C^{\top},\quad 
$$
where $F^{\top}=\begin{bmatrix}f_0 &\ldots & f_{m-1}\end{bmatrix}$, 
$C^{\top}=\begin{bmatrix}c_1 &\ldots & c_{n}\end{bmatrix}$ and where
\begin{equation}
V^{\boldsymbol\ell}_m=\begin{bmatrix} 1 & \alpha_1 & \alpha_1^2 & \ldots & \alpha_1^{m-1}\\
1 & \alpha_2 & \alpha_2^2 & \ldots & \alpha_2^{m-1}\\
\vdots & \vdots & \vdots &&\vdots\\
1 & \alpha_n & \alpha_n^2 & \ldots & \alpha_n^{m-1}\end{bmatrix}=(V_m^{\boldsymbol r})^{\top}
\label{1.6}
\end{equation}
are the left and right Vandermonde matrices associated with the given
$\alpha_1,\ldots,\alpha_n\in\bH$. Recall that the rank of a quaternion matrix is defined as 
the dimension of the left linear span of its rows or equivalently (by \cite[Theorem 7]{lam1}),
as the dimension of the right span of its columns. In general, the ranks of a matrix 
can be different from the rank of its transpose. For a set $\Lambda$, we
will write $\sharp(\Lambda)$ for its cardinality. The following result
is due to T.-Y. Lam \cite{lam1}:

\medskip

\begin{theorem}
Let $S_1,\ldots,S_\ell$ be all distinct conjugacy classes having non-empty
intersection with the set $\Lambda=\{\alpha_1,\ldots,\alpha_n\}$, and let 
\begin{equation}
\kappa=\mu_1+\ldots+\mu_\ell,\quad\mbox{where}\quad \mu_j=\left\{\begin{array}{ccc}
1, & \mbox{if} & \sharp (S_j\cap \Lambda)=1,\\
2,&\mbox{if} & \sharp (S_j\cap \Lambda)\ge 2.\end{array}\right.
\label{1.7}
\end{equation}
Then ${\rm rank} V^{\boldsymbol\ell}_m={\rm rank} V^{\boldsymbol r}_m
={\rm min}(m,\kappa)$. In particular, the square matrix
$V^{\boldsymbol\ell}_n$ is invertible if and only if all elements in $\Lambda$ are distinct 
and none three of them belong to the same conjugacy class.
\label{T:1.1}
\end{theorem}

\smallskip

The result was established in \cite{lam1} in a more general setting of division rings with a fixed
endomorphism and the integer $\kappa$ was identified with the minimally possible degree of a nonzero polynomial
having left zeros at $\alpha_1,\ldots,\alpha_n$. It follows from Theorem \ref{T:1.1} that 
if none three of $\alpha_1,\ldots,\alpha_n$ belong to the same conjugacy class, then
the problems \eqref{1.4} and \eqref{1.5} have (unique) solutions of degree less than $n$ 
for any choice of $c_i$'s. Otherwise, the problems may have no solutions, which is
indicated by the mext result (\cite{genstr}).

\medskip

\begin{theorem}
For $f\in\bH[z]$ and three distinct equivalent elements $\alpha_1,\alpha_2,\alpha_3\in\bH$,
\begin{align}
f^{\bl}(\alpha_3)=&(\alpha_3-\alpha_2)(\alpha_1-\alpha_2)^{-1}f^{\bl}(\alpha_1)
+(\alpha_1-\alpha_3)(\alpha_1-\alpha_2)^{-1}f^{\bl}(\alpha_2),\label{rep1}\\
f^{\br}(\alpha_3)=&f^{\br}(\alpha_1)(\alpha_1-\alpha_2)^{-1}(\alpha_3-\alpha_2)
+f^{\br}(\alpha_2)(\alpha_1-\alpha_2)^{-1}(\alpha_1-\alpha_3).\label{rep2}
\end{align}
\label{T:1.6}
\end{theorem}
Indeed, in order the problem \eqref{1.4} or \eqref{1.5} 
to have a solution, the target values assigned to any triple of equivalent nodes must satisfy 
certain consistency conditions pointed out in Theorem \ref{T:1.6}.
If the problem is consistent, then it admits a (unique) solution of degree less than
${\rm rank} V^{\boldsymbol\ell}_{m}$; 
see \cite{bol} for details.

\bigskip
\noindent
{\bf 1.2. Main results.} Left divided differences of a given polynomial based on the 
ordered collection $\balpha=(\alpha_1,\ldots,\alpha_k)$ of nodes are defined 
in terms of left backward shift operators \eqref{1.3} by 
$[\alpha_1;f]_{\boldsymbol\ell}=f^{\bl}(\alpha_1)$ and 
$$
\left[\alpha_1,\ldots,\alpha_j;f\right]_{\boldsymbol\ell}=
(L_{\alpha_{j-1}}\cdots L_{\alpha_1}f)^{\bl}(\alpha_j)\quad\mbox{for}\quad j\ge 2.
$$
The right divided differences are defined by similar formulas (with right evaluations 
instead of the left and with $R_\alpha$ instead of $L_\alpha$; see \eqref{2.7}).
As we will see in  Section 2, quaternionic divided differences are quite different from their 
complex prototypes, unless the nodes $\alpha_i$'s commute, i.e., $\alpha_i\in\C_{\alpha_1}$ for 
$i=1,\dots,k$. Divided differences based on a spherical chain are of special interest.
Let us introduce notation
\begin{equation}
\Delta_{\boldsymbol\ell}(\balpha;f)=\sbm{[\alpha_1;f]_{\boldsymbol\ell}
\\ [\alpha_1,\alpha_2;f]_{\boldsymbol\ell}\\ \vdots \\
[\alpha_1,\alpha_2,\ldots,\alpha_k;f]_{\boldsymbol\ell}},\quad
{\balpha}=(\alpha_1,\ldots,\alpha_k),\quad f\in\bH[z]
\label{1.10}
\end{equation}
for the column of left divided differences of a given polynomial $f$ based on the spherical
chain $\balpha$. The next theorem shows that the columns
$\Delta_{\boldsymbol\ell}(\boldsymbol\alpha_1;f)$ and
$\Delta_{\boldsymbol\ell}(\boldsymbol\alpha_2;f)$
associated with the chains of the same length $k$ from the same conjugacy class
$S\subset \bH$ and with distinct leftmost entries determine all left divided differences of $f$
of order up to $k$ based on {\em any} chain $\boldsymbol\alpha_3\subset S$.

\medskip

\begin{theorem}
Let $\balpha_i=(\alpha_{i,1},\ldots,\alpha_{i,k})$ ($i=1,2,3$)
be three spherical chains of the same length and in the same conjugacy class $S\subset\bH$ and let
$\alpha_{1,1}\neq\alpha_{2,1}$. Then there exist $k\times k$ matrices $A$ and $B$ such that 
\begin{equation}
\Delta_{\boldsymbol\ell}(\balpha_3;f)=A\Delta_{\boldsymbol\ell}(\balpha_1;f)
+B\Delta_{\boldsymbol\ell}(\balpha_2;f)\quad\mbox{for any}\quad f\in\bH[z].
\label{1.11}
\end{equation}
\label{T:2.4}
\end{theorem}
The proof (along with explicit formulas for $A$ and $B$) will be given in Section 5.

\medskip

\begin{definition}
{\rm Given a spherical chain $\balpha=(\alpha_1,\ldots,\alpha_k)$, we define the 
{\em left confluent Vandermonde matrix} $V^{\boldsymbol\ell}_{m}(\balpha)$ by
\begin{equation} 
V^{\boldsymbol\ell}_{m}(\balpha)=\left[[\alpha_1,\alpha_2,\ldots,\alpha_i; \, z^j]_{\boldsymbol \ell}
\right]_{i=1,\ldots,k}^{j=1,\ldots,m}\in\bH^{k\times m}.
\label{1.12}  
\end{equation}
The left confluent Vandermonde matrix based on $n$ spherical chains
\begin{equation}
\balpha_i=(\alpha_{i,1},\ldots,\alpha_{i,k_i}),\quad i=1,\ldots,n,
\label{1.13}    
\end{equation}  
is given by  
\begin{equation}
V^{\boldsymbol\ell}_{m}(\balpha_1,\ldots,\boldsymbol\alpha_n)=\begin{bmatrix}
V^{\boldsymbol\ell}_{m}(\balpha_1)\\  \vdots \\ V^{\boldsymbol\ell}_{m}(\balpha_n)
\end{bmatrix}   
\label{1.14}    
\end{equation}
where the matrices $V^{\boldsymbol\ell}_{m}(\balpha_j)$ are defined via formula
\eqref{1.12}. Similarly, the right confluent Vandermonde matrix based on spherical chains \eqref{1.13}
is given by
\begin{equation}
V^{\boldsymbol r}_{m}(\balpha_1,\ldots,\balpha_n)=
\begin{bmatrix}
V^{\boldsymbol r}_{m}(\balpha_1)& \ldots & V^{\boldsymbol r}_{m}(\balpha_n)
\end{bmatrix},
\label{1.15}
\end{equation}
where the $m\times k$ block corresponding to a sole chain $\balpha=(\alpha_1,\ldots,\alpha_k)$ is defined 
as
\begin{equation}
V^{\bf r}_{m}(\balpha)=\left[[z^j; \, \alpha_1,\alpha_2,\ldots,\alpha_i]_{\bf r}
\right]_{i=1,\ldots,m}^{j=1,\ldots,k}.
\label{1.16}
\end{equation}}  
\label{D:1.5}
\end{definition}
Vaguely speaking, the confluent Vandermonde matrix should be defined
so that it will be non-singular in situations where the usual Vandermonde matrix is
singular, that is, according to Theorem \ref{T:1.1}, if
$$
\mbox{(1) $\alpha_{i_1}=\alpha_{i_2}\quad$
and/or\quad (2) $\alpha_{i_1}\sim\alpha_{i_2}\sim\alpha_{i_3}$.}
$$
It turns out that 
the matrices \eqref{1.14}, \eqref{1.15} do the job. As we will see in Section 4,
these matrices arise in the context of Lagrange-Hermite interpolation, like regular
Vandermonde matrices do within Lagrange interpolation problems \eqref{1.5}, \eqref{1.6}. 
In case $k_i=1$ for $i=1,\ldots,n$, the matrices \eqref{1.14}, \eqref{1.15} amount 
to $V^{\boldsymbol\ell}_{m}$ and $V^{\bf r}_{m}$ in \eqref{1.6}. If the spherical chain 
$\balpha$ is of the form \eqref{1.9}, then $V^{\boldsymbol\ell}_{m}(\balpha)$
and $V^{\bf r}_{m}(\balpha)$ take the form
\begin{equation}
V^{\boldsymbol\ell}_{m}(\balpha):=\sbm{1&\alpha & \alpha^2 & &\ldots&\ldots&
\alpha^{m-1}\\
0 & 1 & 2\alpha  &&\ldots&\ldots & (m-1)\alpha^{m-2}\\
\vdots & \ddots & \ddots&\ddots &&&\vdots \\
0 &\ldots & 0 & 1 & k\alpha & \ldots&
\frac{(m-1)!}{(n-1)!}\alpha^{m-k}}=V^{\bf r}_{m}(\balpha)^\top.
\label{1.17}
\end{equation}
In the contrast to the basic case \eqref{1.6}, left and right confluent Vandermonde
matrices based on the same noncommutative spherical chains in general are not transposes
of each other (see Example \ref{E:1.3} below). However, they are related as is indicated in the
following result (see Section 2.2 below for the proof).

\medskip

\begin{lemma}
Let $V^{\boldsymbol\ell}_{m}(\balpha_1,\ldots,\balpha_n)$
and $V^{\boldsymbol r}_{m}(\overline{\balpha}_1,\ldots,\overline{\balpha}_n)$
be the left and the right confluent Vandermonde matrices based on on spherical chains
$\balpha_i=(\alpha_{i,1},\ldots,\alpha_{i,k_i})$ and
$\overline{\balpha}_i:=(\overline{\alpha}_{i,1},\ldots,\overline{\alpha}_{i,k_i})$,
respectively. Then
\begin{equation}
V^{\boldsymbol\ell}_{m}(\balpha_1,\ldots,\balpha_n)=
(V^{\boldsymbol r}_{m}(\overline{\balpha}_1,\ldots,\overline{\balpha}_n))^*.
\label{1.18}
\end{equation}
\label{L:1.4}
\end{lemma}
To present the confluent version of Theorem \ref{T:1.1}, we need the analogs of the integers
$\mu_j$'s in \eqref{1.7}. If the conjugacy class $S\subset\bH$ contains only one spherical chain $\balpha_i$
from \eqref{1.13}, we let $\mu(S)=k_i$.
If $S$ contains $d\ge 2$  spherical chains \eqref{1.13}, we pick the longest  chain  
 $\balpha_i=(\alpha_{i,1}\ldots,\alpha_{i,k_i})\subset S$, and for any other chain
$\balpha_j=(\alpha_{j,1}\ldots,\alpha_{j,k_j})\subset S$, define the integer
\begin{equation}
\nu_j=\left\{\begin{array}{ccc}0,&\mbox{if}& \alpha_{j,1}\neq \alpha_{i,1}, \\
\max\{r: \, \alpha_{j,\xi}=\alpha_{i,\xi} \; (1\le \xi\le r)\} , &\mbox{if}&\alpha_{j,1}=\alpha_{i,1}.
\end{array}\right.
\label{1.18a}
\end{equation}
We then let 
\begin{equation}
\mu(S)=k_i+{\displaystyle\max_{j\neq i}\{k_j-\nu_j\}}.
\label{1.18c}
\end{equation}
Note that if there are several chains of the maximal length,
then the values of $\nu_j$ in \eqref{1.18a} depend on which one of the longest chains has
been chosen for the comparison. However, the integer \eqref{1.18c} is independent of this choice
(see Proposition \ref{P:3.4} below).

\medskip

\begin{theorem}
Let $V^{\boldsymbol\ell}_{m}(\balpha_1,\ldots,\balpha_n)$ and 
$V^{\bf r}_{m}(\balpha_1,\ldots,\balpha_n)$ be the confluent 
Vandermonde matrices based on spherical chains \eqref{1.13}. To each conjugacy 
class $S_j$ containing at least one of these chains, assign the integer $\mu(S_j)$ as 
in \eqref{1.18c}, \eqref{1.18a}, and let $\; \kappa=\sum_{j} \mu(S_j)$. Then
\begin{equation}
{\rm rank} \, V_m^{\boldsymbol\ell}(\balpha_1,\ldots,\balpha_n)={\rm min}(m,\kappa)=
{\rm rank} \, V^{\bf r}_{m}(\balpha_1,\ldots,\balpha_n).
\label{1.18b}
\end{equation}
\label{T:1.2}
\end{theorem}
\begin{remark}
{\rm The second equality in \eqref{1.18b} is a consequence of the first, by Lemma \ref{L:1.4}.
Indeed, upon applying formulas \eqref{1.18a} to the chains 
$\overline{\balpha}_1,\ldots,\overline{\balpha}_n$ we come up 
with the same integer $\kappa$ as for the original chains ${\balpha}_1,\ldots,{\balpha}_n$. 
Assuming that the first equality in \eqref{1.18b} holds true,  we conclude that ${\rm rank} \, 
V_m^{\boldsymbol\ell}(\balpha_1,\ldots,\balpha_n)=
{\rm rank} \, V_m^{\boldsymbol\ell}(\overline\balpha_1,\ldots,\overline\balpha_n)$. Since 
${\rm rank} \, A={\rm rank} \, A^*$ for any matrix $A$ over $\bH$, the second equality in \eqref{1.18b} 
follows by \eqref{1.18}.}
\label{R:1.6}
\end{remark}

\smallskip

The outline of the paper is as follows. In Section 2, we present
confluent Vandermonde matrices as unique solutions of certain Stein equations
and prove Lemma \ref{L:1.4}. In Section 3, we recall some basic facts on indecomposable 
polynomials and show that the integer $\kappa$ in Theorem  \ref{T:1.2} can be 
alternatively introduced as the degree of the least right common multiple of polynomials
$P_{\balpha_1},\ldots,P_{\balpha_n}$ associated with the spherical chains $\balpha_1,\
\ldots,\balpha_n$ by formula \eqref{3.4}. In Section 4, we demonstrate how confluent Vandermonde matrices
arise in the context of a Lagrange-Hermite type interpolation problem and present necessary and sufficient 
condition for  a square confluent Vandermonde matrix to be invertible. The proofs of Theorems \ref{T:2.4}
and \ref{T:1.2} are given in Section 5. Several extensions of these results to the setting of 
formal power series over quaternions are presented in the concluding Section 6.

\section{Preliminaries}
\setcounter{equation}{0}

In the complex setting, divided differences 
are defined in terms of the operator
$L_\alpha: \, f(z)\to \frac{f(z)-f(\alpha)}{z-\alpha}$ acting as a backward shift on the
sequence of Taylor coefficients of $f$ at a given point $\alpha\in\C$.
Evaluation formulas for $L_\alpha f$ are:
\begin{equation}
L_\alpha 
f(\beta)=\left\{\begin{array}{ccc}\frac{f(\beta)-f(\alpha)}{\beta-\alpha}&\mbox{if}&
\beta\neq
\alpha,\\
f^\prime(\alpha)&\mbox{if}& \beta= \alpha.\end{array}\right.
\label{2.1}
\end{equation}
Furthermore, the operators $L_\alpha$ and $L_\beta$ commute and satisfy the Hilbert
identity $L_\beta-L_\alpha=(\beta-\alpha)L_\beta L_\alpha$ due to which 
$L_{\alpha_{n-1}}\ldots L_{\alpha_2}L_{\alpha_1}f(\alpha_n)$ is the $n$-th 
divided difference of $f$ at the nodes $\alpha_1,\ldots,\alpha_n$ which in turn, is an ingredient
of the Newton's interpolation formula.

\smallskip

In the noncommutative quaternionic setting we distinguish the left and the right
backward shifts \eqref{1.3}, unless $\alpha$ is real in which case $L_\alpha=R_\alpha$.
If we denote by $f^{(k)}$ the $k$-th formal derivative of
$f\in\bH[z]$, then a straightforward verification shows that for any fixed $\alpha\in\bH$,
\begin{equation}
f=\sum_{k=0}^{\deg (f)} \bp_\alpha^k \frac{(f^{(k)})^{\bl}(\alpha)}{k!}=\sum_{k=0}^{\deg 
(f)}
\frac{(f^{(k)})^{\br}(\alpha)}{k!}\bp_\alpha^k,\quad \bp_\alpha(z):=z-\alpha.
\label{2.2}
\end{equation}
In terms of the latter Taylor expansions, the operators \eqref{1.3}
take the form
$$
L_\alpha f=\sum_{k=0}^{\deg (f)-1} \bp_\alpha^k 
\frac{(f^{(k+1)})^{\bl}(\alpha)}{(k+1)!},\quad
R_\alpha f= \sum_{k=0}^{\deg (f)-1} \frac{(f^{(k+1)})^{\br}(\alpha)}{(k+1)!}\bp_\alpha^k,
$$
which justifies the ``backward shift" terminology. The operators $L_\alpha$ and 
$L_\beta$ are essentially non-commuting.

\medskip

\begin{proposition}
$L_\alpha L_\beta=L_\beta L_\alpha$ if and only if $\alpha\beta=\beta\alpha$, in which case
\begin{equation}
L_\alpha-L_\beta=(\alpha-\beta)L_\alpha L_\beta.
\label{2.3}  
\end{equation}
{\rm Indeed, if $\alpha\beta=\beta\alpha$, then the asserted equalities are verified as in the 
complex case. On the other hand, since 
$\; (L_\alpha L_\beta-L_\beta L_\alpha)(z^4)=\beta\alpha-\alpha\beta$
(by the first formula in \eqref{1.3}), commutation equality $L_\alpha L_\beta-L_\beta L_\alpha$
implies $\beta\alpha=\alpha\beta$. A similar statement holds for right backward shifts}.
\label{P:2.1}
\end{proposition}

\bigskip
\noindent
{\bf 2.1. Divided differences.} Given a polynomial $f\in\bH[z]$, the successive
application of formula \eqref{1.1} to elements $\alpha_1,\ldots,\alpha_{n}\in\bH$
and polynomials $f, \; L_{\alpha_1}f, \; L_{\alpha_2}L_{\alpha_1}f,\ldots$
leads us to the representation
\begin{align}
f=f^{\bl}(\alpha_1)&+\sum_{k=1}^{n-1} \bp_{\alpha_1}\ldots \bp_{\alpha_{k}}\cdot
(L_{\alpha_k}\cdots L_{\alpha_1}f)^{\bl}(\alpha_{k+1})\notag\\
&+\bp_{\alpha_1}\ldots \bp_{\alpha_{n}}\cdot (L_{\alpha_n}\cdots L_{\alpha_1}f),\label{2.4}  
\end{align}
which, being the (left) quaternionic analog of the Newton interpolation formula,
suggests to define {\em quaternionic left divided differences}
\begin{equation}
[\alpha_1;f]_{\boldsymbol\ell}=f^{\bl}(\alpha_1),\qquad 
\left[\alpha_1,\ldots,\alpha_k;f\right]_{\boldsymbol\ell}=
(L_{\alpha_{k-1}}\cdots L_{\alpha_1}f)^{\bl}(\alpha_k)\quad\mbox{for}\quad k\ge 1.
\label{2.5}
\end{equation}
Similarly, the formula based on the successive application of the second representation in
\eqref{1.1},
\begin{align}
f=f^{\br}(\alpha_1)&+\sum_{k=1}^{n-1}
(R_{\alpha_k}\cdots R_{\alpha_1}f)^{\br}(\alpha_{k+1})\cdot \bp_{\alpha_{k}}\ldots
\bp_{\alpha_1}\notag\\
&+(R_{\alpha_n}\cdots R_{\alpha_1}f) \cdot \bp_{\alpha_{n}} \ldots \bp_{\alpha_1}
\label{2.6}
\end{align}
suggests to introduce {\em quaternionic right divided differences} by
\begin{equation}
[f;\alpha_1]_{\boldsymbol r}=f^{\br}(\alpha_1),\qquad
\left[f;\alpha_1,\ldots,\alpha_k\right]_{\boldsymbol r}=
(R_{\alpha_{k-1}}\cdots R_{\alpha_1}f)^{\br}(\alpha_k)\quad\mbox{for}\quad k\ge 1.
\label{2.7}
\end{equation}
Letting $\alpha_j=\alpha$ for $j=1,\ldots,n$ in \eqref{2.4}, \eqref{2.6} and comparing the
obtained representations with \eqref{2.2} we conclude that
\begin{equation}
[\underbrace{\alpha,\ldots,\alpha}_{(k+1) \; {\rm
times}};f]_{\boldsymbol \ell}=\frac{(f^{(k)})^{\bl}(\alpha)}{k!},\quad
[f;\underbrace{\alpha,\ldots,\alpha}_{(k+1) \; {\rm
times}}]_{\boldsymbol r}=\frac{(f^{(k)})^{\br}(\alpha)}{k!}
\label{2.8}
\end{equation}
for all $k\ge 0$. The latter formulas justify equalities \eqref{1.10}.

\smallskip

The difference between the complex and quaternionic settings becomes transparent
even in the case where $k=2$. It is not hard to show that if $\alpha_2\not\sim\alpha_1$, 
then
\begin{equation}
[\alpha_1,\alpha_2;f]_{\boldsymbol \ell}=(\widetilde{\alpha}_2-\alpha_1)^{-1}(f^{\bl}
(\widetilde{\alpha}_2)-f^{\bl}(\alpha_1)), \label{2.9}
\end{equation}
where $\widetilde{\alpha}_2=(\alpha_2-\overline{\alpha}_1)^{-1}\alpha_2
(\alpha_2-\overline{\alpha}_1)$.
The formula \eqref{2.9} is similar to its complex counterpart \eqref{2.1}, but the element 
$\alpha_2$
is replaced by the equivalent element $\widetilde{\alpha}_2$ which is equal to $\alpha_2$ 
if and  only if $\alpha_1$ and $\alpha_2$ commute. If $\alpha_1\sim\alpha_2\neq 
\overline{\alpha}_1$, then (as we will see in Example \ref{E:5.3} below) the formula for 
$[\alpha_1,\alpha_2;f]_{\boldsymbol \ell}$ invokes not only the values of  $f$ at 
$\alpha_1$ and $\alpha_2$, but also the value of $f^\prime$:
\begin{equation}
[\alpha_1,\alpha_2;f]_{\boldsymbol \ell}=(\alpha_2-\overline{\alpha}_2)^{-1}
\left(f^{\bl}(\alpha_2)-f^{\bl}(\alpha_1)+(\alpha_2-\overline{\alpha}_1)
f^{\prime\bl}(\alpha_1)\right)).
\label{2.10}
\end{equation}
Thus, the divided differences based on a spherical chain (different from that in 
\eqref{2.8}) is the object which does not appear in the commutative setting.

\medskip

\begin{lemma}
For any $h\in\bH[z]$, $x\in\R$  and $\alpha_1,\alpha_2,\ldots,\alpha_k\in\bH$,
\begin{equation}
\left[\alpha_1,\ldots,\alpha_k;h\right]_{\boldsymbol \ell}=
(\alpha_k-x)\left[\alpha_1,\ldots,\alpha_k;L_xh\right]_{\boldsymbol
\ell}+\left[\alpha_1,\ldots,\alpha_{k-1};L_xh\right]_{\boldsymbol \ell}.
\label{2.11}
\end{equation}
\label{L:2.2}
\end{lemma}
{\bf Proof:}
Proposition \ref{P:2.1} applies to $\alpha=\alpha_{k-1}$ and $\beta=x$
(since $x$ is real), and the formula \eqref{2.3} implies
$$
L_{\alpha_{k-1}}=(I+(\alpha_{k-1}-x)L_{\alpha_{k-1}})L_x.
$$
Applying formula \eqref{1.1} to $f=L_{\alpha_{k-2}}\cdots L_{\alpha_1}L_xh$
and $\alpha=\alpha_{k-1}$ leads us to
\begin{align*}  
&L_{\alpha_{k-2}}\cdots L_{\alpha_1}L_xh-(z-\alpha_{k-1})L_{\alpha_{k-1}}L_{\alpha_{k-2}}\cdots L_{\alpha_1}L_xh\\
&=(L_{\alpha_{k-2}}\cdots 
L_{\alpha_1}L_xh)^{\bl}(\alpha_{k-1})=\left[\alpha_1,\ldots,\alpha_{k-1};L_xh\right]_{\boldsymbol
\ell}.
\end{align*}
Making use of the two last relations and taking into account that $L_x L_{\alpha_j}= L_{\alpha_j}L_x$
(by Proposition \ref{P:2.1}), we have
\begin{align}
&L_{\alpha_{k-1}}L_{\alpha_{k-2}}\cdots L_{\alpha_1}h-(z-x) L_{\alpha_{k-1}}L_{\alpha_{k-2}}\cdots
L_{\alpha_1}L_xh\notag\\
&=(I+(\alpha_{k-1}-x)L_{\alpha_{k-1}})L_xL_{\alpha_{k-2}}\cdots L_{\alpha_1}h\notag\\
&\quad -(z-\alpha_{k-1})L_{\alpha_{k-1}}L_{\alpha_{k-2}}\cdots L_{\alpha_1}L_xh
-(\alpha_{k-1}-x)L_{\alpha_{k-1}}L_{\alpha_{k-2}}\cdots L_{\alpha_1}L_xh\notag\\
&=L_xL_{\alpha_{k-2}}\cdots L_{\alpha_1}h-L_{\alpha_{k-2}}\cdots L_{\alpha_1}L_xh
+\left[\alpha_1,\ldots,\alpha_{k-1};L_xh\right]_{\boldsymbol \ell}\notag\\
&=\left[\alpha_1,\ldots,\alpha_{k-1};L_xh\right]_{\boldsymbol \ell}.\notag
\end{align}
Evaluating the latter equality at $z=\alpha_k$ on the left gives, in view of \eqref{2.5},
$$
\left[\alpha_1,\ldots,\alpha_k;h\right]_{\boldsymbol
\ell}-(\alpha_k-x)\left[\alpha_1,\ldots,\alpha_k;L_xh\right]_{\boldsymbol
\ell}=\left[\alpha_1,\ldots,\alpha_{k-1};L_xh\right]_{\boldsymbol \ell}
$$
which is equivalent to \eqref{2.11}. \qed

\bigskip
\noindent
{\bf 2.2. Explicit formulas.} In this sction we present the explicit formulas 
for matrices $V^{\boldsymbol\ell}_{m}(\balpha)$ and $V^{\bf r}_{m}(\balpha)$  
defined in \eqref{1.12}, \eqref{1.16}. Here we will not assume   
that the collection of nodes $\balpha=(\alpha_1,\ldots,\alpha_k)$ is a spherical chain.
Given $\balpha=(\alpha_1,\ldots,\alpha_k)$, we let
\begin{equation}
\mathcal J_{{\balpha}}=\left[\begin{array}{ccccc} \alpha_1&0&\ldots&&0\\ 1
&\alpha_2&0&&\\ 0&1&\ddots&\ddots&\vdots\\
\vdots&\ddots&\ddots&\ddots&0\\ 0&\ldots &0&1&\alpha_k
\end{array}\right]\quad\mbox{and}\quad
E_k=\left[\begin{array}{c}1 \\ 0 \\ \vdots \\ 0\end{array}\right]
\label{2.12}
\end{equation}
and we let $v_{i,j}$ denote the $(i,j)$-entry of the matrix $V^{\boldsymbol\ell}_{m}(\balpha)$:
$$
V^{\boldsymbol\ell}_{m}(\balpha)=\left[v_{ij}\right]_{i=1,\ldots,k}^{j=1,\ldots,m},\qquad
v_{i,j}:=[\alpha_1,\alpha_2,\ldots,\alpha_i; \, z^j]_{\boldsymbol \ell}.
$$
Since  $\deg (L_\alpha f)=\deg (f)-1$ for each polynomial $f$ of positive degree 
and since the backward shift of a monic polynomial is again monic (or identical
zero), it follows that $V^{\boldsymbol\ell}_{m}(\balpha)$ is upper triangular 
with all diagonal entries equal one. In particular, the leftmost column $V_1$ of $V^{\boldsymbol\ell}_{m}(\balpha)$
equals $E_k$. Applying equality \eqref{2.11} to $x=0$ and $h=z^j$ and taking into account that
$L_0z^j=z^{j-1}$, we get
\begin{equation}
\left[\alpha_1,\ldots,\alpha_k;z^j\right]_{\boldsymbol \ell}=
\alpha_k\left[\alpha_1,\ldots,\alpha_k;z^{j-1}\right]_{\boldsymbol
\ell}+\left[\alpha_1,\ldots,\alpha_{k-1};z^{j-1}\right]_{\boldsymbol \ell},
\label{2.13}
\end{equation}
which implies the recursion
\begin{equation}
v_{i,j}=\alpha_{i}v_{i,j-1}+v_{i-1,j-1}\quad (i,j\ge 2)
\label{2.14}
\end{equation}
for the entries $v_{i,j}$ of $V^{\boldsymbol\ell}_{m}(\balpha)$. Observe that
$$
v_{1,j}=[\alpha_1;z^j]=\alpha_1^j=\alpha_1 v_{1,j-1},
$$
which together with \eqref{2.14} imply that the  
consecutive columns in $V^{\boldsymbol\ell}_{m}(\balpha)$ are related by
$V_j=\mathcal J_{\balpha}V_{j-1}$ $(j=2,\ldots,m$). Since $V_1=E_k$, the
latter recursion gives $V_j=\mathcal J_{\balpha}^{j-1}E_k$. Similarly, using
relation
$$
\left[z^j;\alpha_1,\ldots,\alpha_k\right]_{\boldsymbol r}=
\alpha_k\left[z^{j-1};\alpha_1,\ldots,\alpha_k\right]_{\bf
r}\alpha_k+\left[z^{j-1};\alpha_1,\ldots,\alpha_{k-1}\right]_{\bf r},
$$
the ``right" counter-part of \eqref{2.13}, one can show that the 
$j$-th row in the matrix \eqref{1.16} equals $E_k^\top (\mathcal J_{{\balpha_i}}^{\top})^{j-1}$.
We thus arrive at the following result.

\medskip

\begin{lemma}
Given $\balpha=(\alpha_1,\ldots,\alpha_k)$, let $V^{\boldsymbol \ell}_m(\balpha)$ and $V^{\boldsymbol r}_m(\balpha)$ be 
defined as in \eqref{1.12}, \eqref{1.16} (confluent Vandermonde matrices if $\balpha$ is a spherical chain). Then
\begin{equation}
V^{\boldsymbol\ell}_{m}(\balpha)=\begin{bmatrix}E_k & & \mathcal
J_{{\balpha}}E_k&&\ldots&&\mathcal J_{{\balpha}}^{m-1}E_k\end{bmatrix},\qquad
V^{\boldsymbol r}_{m}(\balpha)=
\begin{bmatrix}E_{k}^{\top}\\  E_{k}^{\top}\mathcal J_{{\balpha}}^{\top}\\ \vdots \\
E_{k}^{\top}(\mathcal J_{{\balpha}}^{\top})^{m-1}\end{bmatrix},\label{2.15}
\end{equation}
where $\mathcal J_{{\balpha}}$ and $E_k$ are given in \eqref{2.12}.
\label{L:2.1}
\end{lemma}

\medskip

\begin{remark}
{\rm If we let $F_m$ be the $m\times m$ lower triangular Jordan block with zero entries on the main
diagonal ($F_m=\left[\delta_{i-1,j}\right]_{i,j=1}^m$, where $\delta_{ij}$ is the Kronecker
symbol), then we can rewrite formulas \eqref{2.15} as
\begin{equation}
V^{\boldsymbol\ell}_{m}(\balpha)=\sum_{j=0}^{m-1}\mathcal 
J_{{\boldsymbol\alpha}}E_kE_m^{\top}(F_m^{\top})^{j},\quad
V^{\bf r}_{m}(\balpha)=\sum_{j=0}^{m-1}F_m^{j}E_mE_k^{\top}(\mathcal J_{{\boldsymbol\alpha}}^{\top})^{j},
\label{2.20a}
\end{equation}
so that $V^{\boldsymbol\ell}_{m}(\boldsymbol\alpha)$ and $V^{\boldsymbol
r}_{m}(\boldsymbol\alpha)$ are unique solutions to the respective Stein equations
\begin{equation}
V^{\boldsymbol\ell}_{m}(\boldsymbol\alpha)-\mathcal
J_{{\boldsymbol\alpha}}V^{\boldsymbol\ell}_{m}(\boldsymbol\alpha)F_m^{\top}=E_kE_m^{\top},\qquad
V^{\boldsymbol r}_{m}(\boldsymbol\alpha)-F_mV^{\boldsymbol r}_{m}(\boldsymbol\alpha)
{\mathcal J}_{{\boldsymbol\alpha}}^{\top}=E_mE_k^{\top}.\label{2.20}
\end{equation}}
\label{R:2.1}
\end{remark}
{\bf Proof:} Indeed, if $X$ is a solution to the Stein equation
$X=\mathcal J_{{\boldsymbol\alpha}}XF_m^{\top}+E_kE_m^{\top}$, we can iterate this equation to represent
$X$ in the form
$$
X=\sum_{j=0}^{\infty}\mathcal
J_{{\boldsymbol\alpha}}E_kE_m^{\top}(F_m^{\top})^{j}
$$
and then to observe that since $(F_m^{\top})^{m}=0$, the formula for $X$ is the same as that for
$V^{\boldsymbol\ell}_{m}(\balpha)$ in \eqref{2.20a}. The fact that $V^{\bf r}_{m}(\balpha)$ is a unique
solution to the second equation in \eqref{2.20} is justified similarly.\qed

\smallskip

\begin{remark}
{\rm Since left evaluation functionals and left backward shift operators
are right linear, the formula \eqref{1.10} defines a right linear  
operator $f\mapsto \Delta_{\boldsymbol\ell}(\balpha;f)$ acting from $\bH[z]$ into $\bH^k$.
Again, we do not assume that $\balpha=\alpha_1,\ldots,\alpha_k)$ is a spherical chain.
By linearity, we have
$$
\Delta_{\boldsymbol\ell}(\balpha;f)=\sum_{j=0}^{m-1}\sbm{[\alpha_1;z^j]_{\boldsymbol\ell} \\
[\alpha_1,\alpha_2;z^j]_{\boldsymbol\ell}\\ \vdots \\
[\alpha_1,\alpha_2,\ldots,\alpha_k;z^j]_{\boldsymbol\ell}} f_j,\quad\mbox{if}\quad
f(z)=\sum_{j=0}^{m-1}z^jf_j,
$$
which on account of definition \eqref{1.12} and formula \eqref{2.15}, can be written as
\begin{equation}
\Delta_{\boldsymbol\ell}(\balpha;f)=\sum_{j=0}^{m-1}\mathcal J_{\balpha}^{j}E_kf_j=
V^{\boldsymbol\ell}_m(\boldsymbol\alpha)\sbm{f_0 \\ \vdots \\ f_{m-1}}.
\label{4.4}
\end{equation}
}\label{R:2.1a}
\end{remark}

\medskip

{\bf Proof of Lemma \ref{L:1.4}:} Due to definitions \eqref{1.14}, \eqref{1.15}, it suffices to verify
\eqref{1.18} for the left and right confluent Vandermonde matrices based on a sole chain, i.e., to
verify  that
\begin{equation}
V^{\boldsymbol\ell}_{m}(\boldsymbol\alpha)=(V^{\boldsymbol
r}_{m}(\overline{\boldsymbol\alpha}))^*\quad\mbox{for}\quad
{\boldsymbol\alpha}=(\alpha_1,\ldots,\alpha_k).
\label{7.0}
\end{equation}
The latter matrices are unique solutions to equations \eqref{2.20}. Replacing
$\boldsymbol\alpha$ by $\overline{\boldsymbol\alpha}$
and taking adjoints in the second equation in \eqref{2.20} we get (since
$({\mathcal J}_{\overline{\boldsymbol\alpha}}^{\top})^*={\mathcal J}_{{\boldsymbol\alpha}}$)
$$
(V^{\boldsymbol r}_{m}(\overline{\boldsymbol\alpha}))^*
-{\mathcal J}_{{\boldsymbol\alpha}} (V^{\boldsymbol r}_{m}(\overline{\boldsymbol\alpha}))^*
F_m^{\top}=E_kE_m^{\top}.
$$
Therefore, $V^{\boldsymbol\ell}_{m}(\boldsymbol\alpha)$ and $(V^{\boldsymbol
r}_{m}(\overline{\boldsymbol\alpha}))^*$ satisfy the same first equation in \eqref{2.20}, which
has a unique solution so that $V^{\boldsymbol\ell}_{m}(\boldsymbol\alpha)=(V^{\boldsymbol
r}_{m}(\overline{\boldsymbol\alpha}))^*$.\qed

\medskip

\begin{example}
{\rm Let us apply formulas \eqref{2.15}
to the spherical chain ${\balpha}=(\alpha,\beta)$ ($\alpha\sim\beta\neq \overline{\alpha}$)
to get
\begin{align}
V^{\boldsymbol \ell}_{m}(\balpha)&=\begin{bmatrix} 1& \alpha & \alpha^2 & \alpha^3&
\ldots &\alpha^{m-1}\\
0 & 1 & \beta+\alpha &  \beta^2+ \beta \alpha+ \alpha^2&
\ldots &{\displaystyle\sum_{j=0}^{m-2}\beta^j\alpha^{m-j-2}}\end{bmatrix},\label{2.19}\\
(V^{\boldsymbol r}_{m}(\balpha))^{\top}&=\begin{bmatrix} 1& \alpha & \alpha^2 &
\alpha^3  &\ldots &\alpha^{m-1}\\ 0 & 1 & \beta+\alpha &  \beta^2+ \alpha\beta+\alpha^2&
\ldots &{\displaystyle\sum_{j=0}^{m-2}\alpha^{m-j-2}}\beta^j
\end{bmatrix}. \notag
\end{align}
Examining the $(2,4)$-entries shows that $V^{\boldsymbol \ell}_{4}(\balpha)\neq V^{\boldsymbol 
r}_{4}(\balpha)^{\top}$ unless $\alpha\beta=\beta\alpha$ which is possible only if $\beta=\alpha$.
On the other hand, it is readily seen that $V^{\boldsymbol \ell}_{4}(\balpha)=
V^{\boldsymbol r}_{4}(\overline{\balpha})^{*}$ (which agrees with Lemma \ref{L:1.4}).}
\label{E:1.3}
\end{example}

\bigskip
\noindent
{\bf 2.3. Similarities and distinctions with the commutative setting.}
In the commutative setting, confluent Vandermonde matrices can be obtained from 
the basic ones by certain limit procedure, sometimes called {\em the confluence of 
one row into another}. We now examine this procedure on $2\times m$ quaternionic left Vandermonde matrices.
According to the commutative recipe,
we start with the left Vandermonde matrix based on the elements $\alpha$ and $\alpha+\varepsilon$
and then pass to the limit as $\varepsilon\to 0$ in the following product:
\begin{equation}
\lim_{\varepsilon\to 0}\begin{bmatrix}1 & 0 \\ -\varepsilon^{-1} & \varepsilon^{-1}\end{bmatrix}
\begin{bmatrix} 1 & \alpha & \alpha^2 & \ldots & \alpha^{m-1}\\
1 & \alpha+\varepsilon & (\alpha+\varepsilon)^2 & \ldots &
(\alpha+\varepsilon)^{m-1}\end{bmatrix}.
\label{7.1}   
\end{equation}
The two leftmost elements in the bottom row are $0$ and $1$.
To compute other entries, we observe that  
\begin{equation}
\varepsilon^{-1}\left((\alpha+\varepsilon)^k-\alpha^k\right)=\sum_{j=0}^{k-1}\varepsilon^{-1}\alpha^j
\varepsilon\alpha^{k-j-1}+\varepsilon^{-1}R(\alpha,\varepsilon)
\label{7.2} 
\end{equation}
for each  $k\ge 2$, where $R(\alpha,\varepsilon)$ is the sum of all words of length $k$ in letters $\alpha$,
$\varepsilon$ and containing at least two letters $\varepsilon$. Then it follows by the triangle inequality
that
$$
|\varepsilon^{-1}R(\alpha,\varepsilon)|\le |\varepsilon|\sum_{j=0}^{k-2}\binom{k-2}{j}
|\alpha|^j|\varepsilon|^{k-j-2}\to 0, \quad\mbox{as} \quad \varepsilon\to 0.
$$
We now conclude from \eqref{7.2} that the limit \eqref{7.1} exists (for $m\ge 3$) if and only
if the limit $\beta={\displaystyle\lim_{\varepsilon\to
0}\varepsilon^{-1}\alpha\varepsilon}$ exists, that is, {\em never}, unless $\alpha\in\R$
(in which case the confluent matrix takes the form \eqref{1.17}).
However, we may adjust \eqref{7.1} by taking restricted limits (i.e., by sending
$\varepsilon$ to the origin along continuous curves) calling any matrix arising in this way
a confluent Vandermonde matrix. If such a restricted limit $\beta$ exists, it necessarily belongs to
the conjugacy class $[\alpha]$ and the corresponding restricted limit \eqref{7.1} amounts to the matrix
$V^{\boldsymbol\ell}_{m}(\boldsymbol\alpha)$ from \eqref{2.19}. To show that any matrix of the form \eqref{2.19}
(that is, with any $\beta\in[\alpha]$) arises in this way, it suffices to represent $\beta$ in the form
$\beta=h\alpha h^{-1}$ for some $h\neq 0$ and evaluate the limit \eqref{7.1} as
$\varepsilon\to 0$ along the line $\ell_h=\{rh: \, r\in\R\}$.

\section{Indecomposable polynomials}
\setcounter{equation}{0}

In this section, we record several basic facts on quaternion polynomials    
needed in the sequel. As the division algorithm
holds in $\bH[z]$ on either side, any (left or right) ideal in $\bH[z]$ is
principal. We will use notation $\langle h\rangle_{\bf r}:=\left\{hq: \; q\in\bH[z]\right\}$
for the right ideal generated by  $h$, and we will write $f\equiv g \, ({\rm mod}_{\bf r}
\, h)$ in case $(f-g)\in \langle h\rangle_{\bf r}$. Similar notations will be used  
for left ideals. Maximal ideals in $\bH[z]$ are generated by linear polynomials
$\bp_\alpha(z):=z-\alpha$ ($\alpha\in\bH$): it follows from \eqref{1.1} that
\begin{equation}f\in \langle \bp_\alpha\rangle_{\bf r} \; \Leftrightarrow \; f^{\bl}(\alpha)=0
\quad\mbox{and}\quad
f\in \langle \bp_\beta\rangle_{\boldsymbol\ell} \; \Leftrightarrow \; f^{\br}(\beta)=0.
\label{3.1}
\end{equation}
In the latter cases we say that $\alpha$ and $\beta$ are respectively, left and right
zeros of $f$. As was shown in \cite{niven}, any (monic) polynomial $p\in\bH[z]$ of degree $\deg
(p)=k\ge 1$ has a left (a right) zero, which along with \eqref{3.1} implies that $f$ can be
factored into the product of linear factors
\begin{equation}
p=\bp_{\alpha_1}\ldots \bp_{\alpha_{k}},\quad \bp_{\alpha_{j}}(z)=z-\alpha_j.
\label{3.2}
\end{equation}
Equivalences \eqref{3.1} are particular cases of more general statements
\begin{equation}
\begin{array}{cl}
f\in \langle \bp_{\alpha_1}\ldots \bp_{\alpha_{n}}\rangle_{\bf r} \; \Leftrightarrow \;
\left[\alpha_1,\ldots,\alpha_j;f\right]_{\boldsymbol\ell}=0,\\ [1mm]
f\in \langle \bp_{\alpha_1}\ldots \bp_{\alpha_{n}}\rangle_{\boldsymbol \ell} \; \Leftrightarrow
\; \left[f;\alpha_1,\ldots,\alpha_j\right]_{\bf r}=0,\end{array}\quad
\mbox{for} \; \;  j=1,\ldots,k,
\label{3.3}
\end{equation}
which in turn, follow from \eqref{2.4}--\eqref{2.7}.

\smallskip

Given polynomials $f,g\in\bH[z]$, their {\em least right common multiple} $h={\bf lrcm}(f,g)$
and their {\em least left common multiple} $\widetilde h={\bf llcm}(f,g)$
are defined as (unique) monic polynomials such that
$$
\langle
h\rangle_{\bf r}=\langle f\rangle_{\bf r}\cap \langle g\rangle_{\bf r},\quad
\langle\widetilde h\rangle_{\boldsymbol\ell}=\langle f\rangle_{\boldsymbol\ell}\cap \langle
g\rangle_{\boldsymbol\ell}.
$$
Following \cite{ore}, let us say that a polynomial $f$ is {\em indecomposable} if it cannot be
represented as the {\bf lrcm} of its proper left (equivalently, as the {\bf llcm} of its proper right) 
divisors. The latter means that the ideal $\langle h\rangle_{\bf r}$ is irreducible in the sense that it is not 
contained into
two distinct proper right ideals in $\bH[z]$. Various characterizations of
indecomposable polynomials are listed below (see e.g., \cite{bol2} for the proof).

\medskip

\begin{theorem}
Let $p\in\bH[z]$ be factored as in \eqref{3.2}. The following are equivalent:
\begin{enumerate}
\item The ideal $\langle p\rangle_{\bf r}$ is irreducible, i.e., $p$ is indecomposable.
\item The ideal $\langle p\rangle_{\boldsymbol\ell}$ is irreducible.
\item $\boldsymbol\alpha=(\alpha_1,\ldots,\alpha_k)$ is a spherical chain.
\item $\alpha_1$ is the only left zero of $p$.
\item $\alpha_k$ is the only right zero of $p$.
\item \eqref{3.2} is a unique factorization of $p$ into the product of linear factors.
\end{enumerate}
\label{T:3.1}
\end{theorem}

\smallskip

The latter theorem establishes a one-to-one correspondence 
\begin{equation}
\balpha=(\alpha_{1},\ldots,\alpha_{k}) \mapsto P_{\balpha}=\bp_{\alpha_{1}}\ldots \bp_{\alpha_{k}}
\label{3.4}
\end{equation}
between spherical chains and  monic indecomposable polynomials. 

\smallskip

Given $f,g\in\bH[z]$, their {\em greatest left common divisor} ${\bf glcd}(f,g)$
is defined as a monic polynomial $d$ of the highest
possible degree such that $f=d \widetilde{f}$ and $g=d \widetilde{g}$ for some $\widetilde{f}, 
\widetilde{g}\in\bH[z]$ or equivalently, as a monic generator of
the right ideal $\langle f\rangle_{\bf r}+\langle g\rangle_{\bf r}$.

\medskip

\begin{remark}
{\rm It follows from the equivalence $(3)\Leftrightarrow(6)$ in Theorem \ref{T:3.1} that
given two spherical chains $\balpha_i=(\alpha_{i,1},\ldots,\alpha_{i,k_i})$ ($i=1,2$)
belonging to the same conjugacy class and sharing $\nu$ leftmost elements, the 
{\bf glcd} of the associated polynomials $P_{\balpha_1}$ and $P_{\balpha_2}$ (see \eqref{3.4}) equals 
$\bp_{\alpha_1}\cdots \bp_{\alpha_\nu}$.}
\label{R:3.2}
\end{remark}

\smallskip

Let us recall that the {\em characteristic polynomial} of a nontrivial conjugacy class $S\subset\bH$
is defined by
\begin{equation}
\cX_S(z)=(z-\alpha)(z-\overline{\alpha})=z^2-z(\alpha+\overline{\alpha})+|\alpha|^2,
\label{3.5}
\end{equation}
where $\alpha$ is any element in $S$; it follows from characterization \eqref{1.0} that formula \eqref{3.5} 
does not depend on the choice of $\alpha\in S$. Since $\cX_S$ is the polynomial of the minimally possible degree 
such that its zero set (left and right, as $\cX\in{\mathbb R}[z]$) coincides with $S$, it is also called
the {\em minimal polynomial} of $S$. We now recall a result from \cite{bol2} concerning least common 
multiples of indecomposable polynomials having zeros in the same conjugacy class. 

\medskip

\begin{lemma}
Given indecomposable polynomials
\begin{equation}
P_{\balpha_i}=\bp_{\alpha_{i,1}}\ldots \bp_{\alpha_{i,k_i}},\qquad
\balpha_i=(\alpha_{i,1},\ldots,\alpha_{i,k_i}),
\label{3.6}
\end{equation}
based on spherical chains $\balpha_1,\ldots\balpha_d$ in the same conjugacy class $S$
such that $\deg(P_{\balpha_1})\ge \deg(P_{\balpha_2})\ge \ldots \ge \deg(P_{\balpha_d})$
(i.e., $k_1\ge k_2\ge \ldots \ge k_d$), let
\begin{equation}
P_{\balpha_j}=p_jh_j,\quad\mbox{where}\quad p_j={\bf glcd}(P_{\balpha_j},P_{\balpha_1})\quad\mbox{for}\quad
j=2,\ldots,d,
\label{3.7}
\end{equation}
and let $m={\displaystyle\max_{2\le j\le d} \deg (h_j)}$. Then
\begin{equation}
{\bf lrcm}(P_{\balpha_1},\ldots,P_{\balpha_d})=
\left\{\begin{array}{ccc} \cX_S^m , &\mbox{if} & m=k_1,\\
\cX_S^m \bp_{\alpha_{1,1}}\bp_{\alpha_{1,2}}\cdots
\bp_{\alpha_{k_1-m}}, &\mbox{if} & m<k_1.\end{array}\right.
\label{3.8}  
\end{equation} 
\label{L:3.3}
\end{lemma}
\begin{corollary}
If $\balpha_{i}=(\alpha_{i,1},\ldots,\alpha_{i,k_i})$ ($i=1,2$) are two spherical chains in the same conjugacy 
class and if $\alpha_{1,1}\neq \alpha_{2,1}$, then 
$\deg({\bf lrcm}(P_{\balpha_1},P_{\balpha_2}))=
\deg(P_{\balpha_1})+\deg(P_{\balpha_2})$.
\label{C:3.5}
\end{corollary}

\smallskip

Since all polynomials \eqref{3.6} are indecomposable, it follows by Remark \ref{R:3.2} that 
$\deg (p_j)=\nu_j$, where $\nu_j$ is the integer defined in \eqref{1.18a} (with $i=1$ since 
in Lemma \ref{L:3.3} we assumed that $\balpha_1$ is the longest chain). Therefore, the integer 
$m$ in \eqref{3.8} equals ${\displaystyle\max_{2\le j\le d}(k_j-\nu_j)}$.
By formula \eqref{3.8}, 
$$
\deg({\bf lrcm}(P_{\balpha_1},\ldots,P_{\balpha_d}))=k_1+m.
$$
Comparing the latter equality with \eqref{1.18c} leads us to the following 

\medskip

\begin{proposition}
{\rm The integer $\mu(S)$ defined in \eqref{1.18c}, \eqref{1.18a} is equal to
$\deg({\bf lrcm}(P_{\balpha_i}: \; \balpha_i\subset S))$}.
\label{P:3.4}
\end{proposition}

\smallskip

We next observe that if $f_1,\ldots,f_n\in\bH[z]$ are ``spherically coprime" (i.e., no zeros 
of $f_i$ and $f_j$ belong to the same conjugacy class), then 
$\deg ({\bf lrcm}(f_1,\ldots,f_n))=\sum_{j=1}^n \deg (f_j)$. 
Combining this observation with respectively, Corollary \ref{C:3.5} and Proposition \ref{P:3.4} leads
us to the following conclusions.

\medskip

\begin{remark}
{\rm If the leftmost elements $\alpha_{1,1},\ldots,\alpha_{n,1}$ in the chains \eqref{1.13} are all
distinct and none three of them belong to the same conjugacy class, then 
$$
\deg({\bf lrcm}(P_{\balpha_1},\ldots,P_{\balpha_n}))=\sum_{j=1}^n \deg (P_{\balpha_j})=\sum_{j=1}^n 
k_j.
$$}
\label{R:3.6}
\end{remark}
\begin{remark}
{\rm The integer $\kappa$ in Theorem \ref{T:1.2} equals
$\deg({\bf lrcm}(P_{\balpha_1},\ldots,P_{\balpha_n}))$.}
\label{R:3.7}
\end{remark}

\section{Lagrange-Hermite interpolation}
\setcounter{equation}{0}
In Section 1, we recalled how Vandermonde matrices arise in the context of
the Lagrange interpolation problem. In this section we formulate and solve
a more general Lagrange-Hermite problem and will see (in formulas \eqref{4.5} and \eqref{4.7} below)
how confluent Vandermonde matrices arise in this more general setting.

\smallskip

Due to equivalence \eqref{3.1}, the typical interpolation condition $f^{\bl}(\alpha)=c$ 
in the left Lagrange problem \eqref{1.4} can be written in the form $f\equiv c \, ({\rm mod}_{\bf r} \, 
\bp_{\alpha})$ where the target value $c$ is understood as an element of $\bH[z]$ of degree zero. 
This form suggests to
consider a more general condition 
\begin{equation}
f\equiv h \, ({\rm mod}_{\bf r} \, p)
\label{4.1} 
\end{equation}
where $p$ is a given indecomposable monic polynomial and $h$ is a
given polynomial of degree $\deg \, (h)<\deg \, (p)=k$. By Proposition \ref{P:2.1}
we may take $p$ in the form \eqref{3.2} where $\boldsymbol\alpha=(\alpha_{1},\ldots,\alpha_{k})$
is a spherical chain, and we may apply representation \eqref{2.4} and notation \eqref{2.5} to write
$$
h=h_{1}+\sum_{j=1}^{k-1}\bp_{\alpha_{1}}\ldots \bp_{\alpha_{j}}\cdot h_{j+1},
\quad\mbox{where}\quad h_{j}=[\alpha_{1},\ldots,\alpha_{j};h]_{\boldsymbol\ell}.
$$
Since any $f\in\bH[z]$ can be represented as
$$
f=[\alpha_{1};f]_{\boldsymbol\ell}+\sum_{j=1}^{k-1}\bp_{\alpha_{1}}\ldots \bp_{\alpha_{j}}\cdot
[\alpha_{1},\ldots,\alpha_{j};f]_{\boldsymbol\ell}+pq\quad \mbox{for some} \; q\in\bH[z],
$$
we conclude that condition \eqref{4.1} can be equivalently written as
\begin{equation}
[\alpha_{1},\ldots,\alpha_{j};f]_{\boldsymbol\ell}=h_{j} \quad\mbox{for}\quad j=1,\ldots,k,
\label{4.2}
\end{equation}
or in the vector form, upon making use of notation \eqref{1.10}, as
\begin{equation}
\Delta_{\boldsymbol\ell}(\boldsymbol\alpha;f)={\rm Col}_{1\le j\le k}h_j.
\label{4.3}
\end{equation} 
Upon invoking formula \eqref{4.4} for $\Delta_{\boldsymbol\ell}(\boldsymbol\alpha;f)$, we 
write interpolation condition \eqref{4.3} in terms of coefficients of the unknown $f$ as 
\begin{equation}
V^{\boldsymbol\ell}_m(\boldsymbol\alpha)\sbm{f_0 \\ \vdots \\
f_{m-1}}=\sbm{h_1 \\ \vdots \\ h_k}.
\label{4.5}
\end{equation}
Similarly, given an indecomposable polynomial $\widetilde{p}=\bp_{\alpha_k}\ldots \bp_{\alpha_1}$
(observe that according to definition \eqref{1.8}, if $\balpha=(\alpha_1,\ldots,\alpha_k)$
is a spherical chain, then the reversed tuple $(\alpha_k,\ldots,\alpha_1)$ is also a spherical chain)
and given $\widetilde{h}\in\bH[z]$ ($\deg (\widetilde{h})<k$), the 
``right" analog of condition \eqref{4.1} is equivalent to the
matrix equation
$$
f\equiv \widetilde{h} \, ({\rm mod}_{\boldsymbol\ell}\,  \widetilde{p})\; \Longleftrightarrow \; 
\begin{bmatrix}f_0 & \ldots & f_{m-1}\end{bmatrix}
V^{\bf r}_m(\boldsymbol\alpha)=\begin{bmatrix}\widetilde{h}_1 & \ldots & 
\widetilde{h}_k\end{bmatrix}
$$
where $\widetilde{h}_j=[\widetilde{h};\alpha_1,\ldots,\alpha_j]_{\bf r}$ for $j=1,\ldots,k$.

\medskip

\begin{remark}
{\rm For the spherical chain \eqref{1.9}, interpolation conditions \eqref{4.2}
prescribe the values of $(f^{(j)})^{\bl}(\alpha)$ for $j=0,\ldots,k-1$ (see formulas \eqref{2.8})
bringing up therefore, the literal quaternionic analog of the well-known Lagrange-Hermite interpolation problem.
In the general case, as is indicated in \eqref{2.10}, conditions \eqref{4.2} prescribe certain
combinations of the values of $f$ and its derivatives at the elements of the chain $\balpha$.
Since the ideal $\langle p\rangle_{\bf r}$ is irreducible, the values of $f$ and its derivatives
cannot be separated, so we indeed have an interpolation problem which does not appear in the 
commutative case. }
\label{R:4.1}
\end{remark}

\medskip
\noindent
{\bf Left Lagrange-Hermite interpolation problem:}  
{\em given $n$ spherical chains \eqref{1.13} and given elements $c_{ij}\in\bH$, find a polynomial 
$f\in\bH[z]$ such that 
\begin{equation}
[\alpha_{i,1},\ldots,\alpha_{i,j};f]_{\boldsymbol\ell}=c_{i,j} \quad\mbox{for}\quad 
i=1,\ldots,n; \; j=1,\ldots,k_i.
\label{4.6}
\end{equation}}
Similarly to \eqref{4.5}, we can write conditions \eqref{4.6} in the matrix form as 
\begin{equation}
V^{\boldsymbol\ell}_m(\balpha_1,\ldots,\balpha_n)\sbm{f_0 \\ \vdots \\
f_{m-1}}=\sbm{C_1 \\ \vdots \\ C_n},
\qquad C_i=\sbm{c_{i,1} \\ \vdots \\ c_{i,k_i}} \; \; (i=1,\ldots,n),   
\label{4.7}
\end{equation}
where $V^{\boldsymbol\ell}_m(\balpha_1,\ldots,\balpha_n)$ is defined in \eqref{1.10}. 

\medskip

\begin{lemma}
Given spherical chains \eqref{1.13} and associated polynomials \eqref{3.6},
the following are equivalent:
\begin{enumerate}
\item $V^{\boldsymbol\ell}_m(\balpha_1,\ldots,\balpha_n)X_m=0, \; $ where 
$X_m=\begin{bmatrix}g_0 & \ldots & g_{m-1}\end{bmatrix}^\top$.
\item The polynomial $g(z)=\sum_{j=1}^{m-1}z^jg_j$ satisfies homogeneous interpolation conditions
\begin{equation}
[\alpha_{i,1},\ldots,\alpha_{i,j};\, g]_{\boldsymbol\ell}=0 \quad\mbox{for}\quad
i=1,\ldots,n; \; j=1,\ldots,k_i.
\label{4.10}  
\end{equation}
\item $g$ is in the right ideal $\langle G\rangle_{\bf r}$ generated by $G={\bf 
lrcm}(P_{\balpha_1},\ldots,P_{\balpha_n})$.
\end{enumerate}
\label{L:4.2}
\end{lemma}

\smallskip
\noindent
{\bf Proof:}
The equivalence $(1)\Leftrightarrow(2)$ is the homogeneous version of the equivalence 
\eqref{4.6}$\Leftrightarrow$\eqref{4.7} above. Furthermore, according to \eqref{3.3}, conditions 
\eqref{4.10} are equivalent to 
$$
g\equiv 0 \, ({\rm mod}_{\bf r} \, P_{\balpha_i})\quad\mbox{for}\quad i=1,\ldots,n,
$$
where the polynomials $P_{\balpha_i}$ are defined in \eqref{3.6}. The latter conditions are
in turn equivalent to the single condition  $g\equiv 0 \, ({\rm mod}_{\bf r} \, G)$. \qed

\medskip

\begin{corollary}
Given spherical chains \eqref{1.13}, let $\kappa$ be the integer defined in Theorem \ref{T:1.2}
(equivalently, $\kappa=\deg({\bf lrcm}(P_{\balpha_1},\ldots,P_{\balpha_n}))$). Then 
$\kappa$ leftmost columns in  
$V^{\boldsymbol\ell}_m(\balpha_1,\ldots,\balpha_n)$ ($m\ge \kappa$) are right linearly independent.
\label{C:4.3}
\end{corollary}

\smallskip
\noindent
{\bf Proof:}
Since $g\equiv 0$ is the only polynomial in the ideal 
$\langle{\bf lrcm}(P_{\balpha_1},\ldots,P_{\balpha_n})\rangle_{\bf r}$ of degree less than $\kappa$, 
it follows by the 
equivalence $(1)\Leftrightarrow(3)$ in Lemma \ref{L:4.2} that the homogeneous equation 
$V^{\boldsymbol\ell}_\kappa(\balpha_1,\ldots,\balpha_n)X_\kappa=0$ has only trivial solution.
Therefore the columns in the matrix $V^{\boldsymbol\ell}_\kappa(\balpha_1,\ldots,\balpha_n)$
(i.e., $\kappa$ leftmost columns in
$V^{\boldsymbol\ell}_m(\balpha_1,\ldots,\balpha_n)$ for $m\ge \kappa$) are right linearly independent.
\qed

\medskip

\begin{theorem}
Assume that the leftmost elements $\alpha_{1,1},\ldots,\alpha_{n,1}$ in given spherical chains \eqref{1.13}
are all distinct and none three of them belong to the same conjugacy class.
Then  the square matrix $V^{\boldsymbol\ell}_N(\balpha_1,\ldots,\balpha_n)$ ($N=\sum_ik_i$)
is invertible and all solutions to the problem \eqref{4.6} are given by the formula
\begin{equation}
f=\widetilde f+Gh,\quad\mbox{where}\quad G:={\bf lrcm}\{ P_{\balpha_1},\ldots,P_{\balpha_n}\},
\label{4.8}
\end{equation}
where $P_{\balpha_i}$ is the polynomial defined in \eqref{3.6}, $h\in\bH[z]$ is a free parameter, and where
\begin{equation}
\widetilde f(z)=\begin{bmatrix}1 & z & \ldots &
z^{N-1}\end{bmatrix}V^{\boldsymbol\ell}_N(\balpha_1,\ldots,\balpha_n)^{-1}\sbm{C_1 \\
\vdots \\ C_n}
\label{4.9}
\end{equation}
is a (unique) solution to the problem \eqref{4.6} of degree less than $N$.
\label{T:4.1}
\end{theorem}

\smallskip
\noindent
{\bf Proof:}
By Remark \ref{R:3.6}, $\deg (G)=N$ and hence, the columns of the matrix 
$V^{\boldsymbol\ell}_N(\balpha_1,\ldots,\balpha_n)$ are right
linearly independent (by Corollary \ref{C:4.3}), so that the matrix is invertible.
We now can solve the non-homogeneous equation \eqref{4.7} (with $m=N$) to get the column of 
coefficients $f_0,\ldots,f_N$ which being substituted into the formula
$$
\widetilde{f}(z)=\sum_{j=0}^{N-1}z^jf_j=\begin{bmatrix}1 & \ldots &
z^{N-1}\end{bmatrix}\sbm{f_0 \\ \vdots \\ f_{N-1}},
$$
leads us to \eqref{4.9}. Since $Gh$ is the general solution to the homogeneous
problem \eqref{4.10}, the formula \eqref{4.8} follows.
Since $\deg (\widetilde{f})<N$ (by \eqref{4.9}) and $\deg (G)=N$,
it follows that  $\deg (\widetilde{f}+Gh) =\deg (G) +\deg(h)\ge N$ for any $h\not\equiv 0$, and therefore, 
$\widetilde{f}$ is indeed a unique solution to the problem \eqref{4.6} of  degree less than $N$. \qed

\medskip

\begin{corollary}
The square confluent Vandermonde matrix $V^{\boldsymbol\ell}_m(\balpha_1,\ldots,\balpha_n)$
based on $n$ spherical chains \eqref{1.13} is invertible if and
only if the leftmost elements $\alpha_{1,1},\alpha_{2,1},\ldots,\alpha_{n,1}$
are all distinct and none three of them belong to the same conjugacy class.
\label{C:4.2}
\end{corollary}
\smallskip
\noindent
{\bf Proof:}
The ``if" part has been proven in Theorem \ref{T:4.1}. To prove the
``only if" part, let us consider the usual Vandermonde matrix
$V_{n}^{\boldsymbol\ell}=\left[\alpha_{i,1}^{k-1}\right]_{i,j=1}^n$
based on the leftmost elements $\alpha_{1,1},\alpha_{2,1},\ldots,\alpha_{n,1}$
of the given spherical chains. If the confluent matrix
$V^{\boldsymbol\ell}_m(\balpha_1,\ldots,\balpha_n)$ is invertible, the rows
in $V_m^{\boldsymbol\ell}$ are left linearly independent, so that  $V_m^{\boldsymbol\ell}$
is invertible as well. Then the claim follows from Theorem \ref{T:1.1}.\qed

\smallskip

The formulation of the right Lagrange-Hermite problem and the ``right" versions of Theorem \ref{T:4.1} and 
Corollary \ref{C:4.2} will be omitted.

\section{Representation formulas for divided differences} 
\setcounter{equation}{0}
Formulas \eqref{rep1}, \eqref{rep2} (termed ``representation formulas" in \cite{genstr}) 
express the value of a polynomial at any point in the conjugacy class
in terms of its values at any two other points from the same class. 
It is clear that the same formulas hold for formal derivatives of $f$, that is,
for divided differences of $f$ based on the spherical chains of the special form \eqref{2.8}. 
Theorem \ref{2.4} asserts that the divided differences of $f$ based on two spherical chains 
of the same length $k$ from the same conjugacy class $S\subset \bH$ and with distinct leftmost entries 
define all divided differences of $f$ of order up to $k$ based on {\em any} chain 
$\boldsymbol\alpha_3\subset S$. A more general question is: {\em given arbitrary spherical chains 
$\balpha_1,\balpha_2\subset S$, for which spherical chains $\boldsymbol\alpha_3\subset S$ does the 
formula \eqref{1.1} hold with some matrices $A$ and $B$ independent of $f$?}

\medskip

\begin{theorem}
Let $\balpha_i=(\alpha_{i,1},\ldots,\alpha_{i,k_i})$ ($i=1,2$, $k_2\le k_1$) be two spherical chains in 
the same
conjugacy class  $S\subset\bH$ sharing $\nu_2$ leftmost elements:
\begin{equation}
\alpha_{1,j}=\alpha_{2,j}\quad(1\le j\le \nu_2)\quad\mbox{and}\quad
\alpha_{1,\nu_2+1}\neq \alpha_{2,\nu_2+1}.
\label{6.3}
\end{equation}
Let $\balpha_3=(\alpha_{3,1},\ldots,\alpha_{3,k_3})$ be another spherical chain in $S$ such that
\begin{equation}
k_3\le k_1\quad\mbox{and}\quad k_3-\nu_3 \le k_2-\nu_2,
\label{6.4}
\end{equation}
where $\nu_3$ is the number of leftmost elements shared by $\balpha_1$ and
$\balpha_3$. Then there exist matrices $A\in\bH^{k_3\times k_1}$ and $B\in\bH^{k_3\times k_2}$ such 
that
\begin{equation}
\Delta_{\boldsymbol\ell}(\balpha_3; \, g)=A\Delta_{\boldsymbol\ell}(\balpha_1; \, g)+
B\Delta_{\boldsymbol\ell}(\balpha_2; \, g)
\label{6.5}
\end{equation}
for any polynomial $g\in\bH[z]$. Consequently,
\begin{equation}
V_m^{\boldsymbol\ell}(\balpha_3)=AV_m^{\boldsymbol\ell}(\balpha_1)+
BV_m^{\boldsymbol\ell}(\balpha_2)\quad \mbox{for}\quad m=1,,2,\ldots.
\label{6.6}
\end{equation}
\label{T:6.2}
\end{theorem}

\smallskip
\noindent
{\bf Proof:}
 We start with the Lagrange-Hermite interpolation problem
\begin{equation}
\begin{array}{ll}
[\alpha_{1,1},\ldots,\alpha_{1,j}; \, f]_{\boldsymbol\ell}&=c_{1,j}\quad\mbox{for}\quad 
j=1,\ldots,k_1,\\ [1mm]
[\alpha_{2,1},\ldots,\alpha_{2,j}; \, f]_{\boldsymbol\ell}&=c_{2,j}\quad\mbox{for}\quad  
j=\nu_2+1,\ldots,k_2
\end{array}
\label{6.7}
\end{equation}
based on the chains $\balpha_1$, $\balpha_2$; the target values $c_{i,j}$ will be specified later.

\medskip
\noindent
{\bf Step 1:} {\em The polynomial $f(z)={\displaystyle\sum_{j=0}^{m-1}z^jf_j}$ satisfies conditions \eqref{6.7} if and 
only if
\begin{equation}
\begin{bmatrix}I_{k_1} & 0 & 0 \\ 0 & 0& I_{k_2-\nu_2}\end{bmatrix}
V_{m}^{\boldsymbol\ell}(\balpha_1,\balpha_2)X_m=\begin{bmatrix}C_1 \\ C_2\end{bmatrix},
\label{6.8}
\end{equation}
where $V_{m}^{\boldsymbol\ell}(\balpha_1,\balpha_2)$ is the confluent Vandermonde matrix based on 
$\balpha_1$, $\balpha_2$,
$$
C_1=\begin{bmatrix}c_{1,1}\\ \vdots \\ c_{1,k_1}\end{bmatrix},\quad
C_2=\begin{bmatrix}c_{2,\nu_2+1}\\ \vdots \\ c_{2,k_2}\end{bmatrix}\quad\mbox{and}\quad
X_m=\begin{bmatrix}f_0 \\ \vdots \\ f_{m-1}\end{bmatrix}.
$$}
Indeed, making use of notation \eqref{1.10}, we write conditions \eqref{6.7} in the vector form as
$$
\Delta_{\boldsymbol\ell}(\balpha_1; \, f)=C_1,\quad
\begin{bmatrix}0 & I_{k_2-\nu_2}\end{bmatrix}\Delta_{\boldsymbol\ell}(\balpha_2; \,
f)=C_2,
$$
or in the matrix form, as
\begin{equation}
\begin{bmatrix}I_{k_1} & 0 & 0 \\ 0 & 0& I_{k_2-\nu_2}\end{bmatrix}
\begin{bmatrix}\Delta_{\boldsymbol\ell}(\balpha_1; \, f)\\ \Delta_{\boldsymbol\ell}(\balpha_2; \,
f)\end{bmatrix}=\begin{bmatrix}C_1 \\ C_2\end{bmatrix},
\label{6.9}
\end{equation}  
and therefore, in the form \eqref{6.8}, due to \eqref{4.4}.

\medskip
\noindent
{\bf Step 2:} {\em The square matrix
\begin{equation}
K=\begin{bmatrix}I_{k_1} & 0 & 0 \\ 0 & 0& I_{k_2-\nu_2}\end{bmatrix}
V_{k_1+k_2-\nu_2}^{\boldsymbol\ell}(\balpha_1,\balpha_2)\quad\mbox{is invertible}.
\label{6.10}
\end{equation}}
If we set $c_{i,j}=0$ everywhere in \eqref{6.7},
the obtained homogeneous conditions are equivalent (due to assumptions \eqref{6.3} and the equivalence
\eqref{3.3}) to
$$
f\equiv 0 ({\rm mod}_{\bf r}P_{\balpha_i})\quad\mbox{where}\quad
P_{\balpha_i}=\bp_{\alpha_{i,1}}\ldots
\bp_{\alpha_{i,k_i}} \; \; (i=1,2).
$$
Equivalently, $f\equiv 0 ({\rm mod}_{\bf r} G)$ where $G={\bf
lrcm}(P_{\balpha_1},P_{\balpha_2})$, and it follows from \eqref{6.3} by Lemma \ref{L:3.3} that
\begin{equation} 
G=\cX_S^{k_2-\nu_2}\bp_{\alpha_{1,1}}\cdots\bp_{\alpha_{1,k_1-k_2+\nu_2}}.
\label{6.11}
\end{equation}
Since $\deg (\cX_S)=2$, we have $\deg (G)=k_1+k_2-\nu_2$ and hence $f\equiv 0$ is the only solution to 
the homogeneous problem \eqref{6.7} of degree less than $k_1+k_2-\nu_2$. Equivalently, the 
homogeneous
matrix equation
$$
\begin{bmatrix}I_{k_1} & 0 & 0 \\ 0 & 0& I_{k_2-\nu_2}\end{bmatrix}
V_{k_1+k_2-\nu_2}^{\boldsymbol\ell}(\balpha_1,\balpha_2)X_{k_1+k_2-\nu_2}=0
$$
has only trivial solution. Therefore, the matrix $K$ of this system
is invertible.

\medskip   
\noindent
{\bf Step 3:} {\em If $f\in\bH[z]$ satisfies homogeneous conditions \eqref{6.7}, then
it also satisfies $\Delta_{\boldsymbol\ell}(\balpha_3; \, f)=0$
for any spherical chain $\balpha_3$ subject to \eqref{6.4}.}

\smallskip

Due to assumptions \eqref{6.4}, the polynomial $P_{\balpha_3}=\bp_{\alpha_{3,1}}\cdots 
\bp_{\alpha_{3,k_3}}$
is a left divisor of the polynomial $G$ given in \eqref{6.11}. Indeed, since $\cX_S$ has real coefficients, it 
commutes with any polynomial in $\bH[z]$. On the other hand, $\cX_S=\bp_\alpha\bp_{\overline{\alpha}}$ 
for each $\alpha\in S$. Since $k_1-k_2+\nu_2\ge k_3-k_2+\nu_2\ge \nu_3$, we have
\begin{align}
G=&\cX_S^{k_2-\nu_2}\bp_{\alpha_{1,1}}\cdots\bp_{\alpha_{1,k_1-k_2+\nu_2}}\notag \\
=&\bp_{\alpha_{1,1}}\cdots\bp_{\alpha_{1,\nu_3}}\cX_S^{k_2-\nu_2}\bp_{\alpha_{1,\nu_3+1}}\cdots 
\bp_{\alpha_{1,k_1-k_2+\nu_2}}
\notag\\
=&\bp_{\alpha_{3,1}}\cdots\bp_{\alpha_{3,\nu_3}}
\bp_{\alpha_{3,\nu_3+1}}\cdots\bp_{\alpha_{3,k_3}}\bp_{\overline{\alpha}_{3,k_3}}
\cdots\bp_{\overline{\alpha}_{3,\nu_3+1}}\times \notag\\
&\qquad \times \cX_S^{k_2-\nu_2-k_3+\nu_2}\bp_{\alpha_{1,\nu_3+1}}\cdots 
\bp_{\alpha_{1,k_1-k_2+\nu_2}}\notag\\
=&P_{\balpha_3} \cX_S^{k_2-\nu_2-k_3+\nu_2}\bp_{\overline{\alpha}_{3,k_3}}
\cdots\bp_{\overline{\alpha}_{3,\nu_3+1}}\bp_{\alpha_{1,\nu_3+1}}\cdots \bp_{\alpha_{1,k_1-k_2+\nu_2}}.
\label{6.12}
\end{align}
If $f\in\bH[z]$ satisfies homogeneous conditions \eqref{6.7}, then $f$ belongs to $\langle 
G\rangle_{\bf r}$,
by (the proof of) Step 2. Then $f\in\langle P_{\balpha_3}\rangle_{\bf r} $ (by \eqref{6.12}) which is
equivalent to equality $\Delta_{\boldsymbol\ell}(\balpha_3; \, f)=0$ (by \eqref{3.2}).

\medskip

We now complete the proof of the theorem.  Given a polynomial $g\in\bH[z]$, we let
\begin{equation}
\begin{array}{ll}
c_{1,j}&=[\alpha_{1,1},\ldots,\alpha_{1,j}; \, g]_{\boldsymbol\ell}\quad\mbox{for}\quad 
j=1,\ldots,k_1,\\ [1mm]
c_{2,j}&=[\alpha_{2,1},\ldots,\alpha_{2,j}; \, g]_{\boldsymbol\ell}\quad\mbox{for}\quad  
j=\nu_2+1,\ldots,k_2,
\end{array}
\label{6.13}
\end{equation}
and consider the interpolation problem \eqref{6.7} based on this data. By Step 1, the problem is 
equivalent to the matrix equation \eqref{6.8} with
$$
\begin{bmatrix}C_1 \\ C_2\end{bmatrix}=
\begin{bmatrix}I_{k_1} & 0 & 0 \\ 0 & 0& I_{k_2-\nu_2}\end{bmatrix}
\begin{bmatrix}\Delta_{\boldsymbol\ell}(\balpha_1; \, g)\\ \Delta_{\boldsymbol\ell}(\balpha_2; \,
g)\end{bmatrix}.
$$
For $m=k_1+k_2-\nu_2$, this equation takes the form
$$
KX_{k_1+k_2-\nu_2}=\begin{bmatrix}I_{k_1} & 0 & 0 \\ 0 & 0& I_{k_2-\nu_2}\end{bmatrix}
\begin{bmatrix}\Delta_{\boldsymbol\ell}(\balpha_1; \, g)\\ \Delta_{\boldsymbol\ell}(\balpha_2; \,
g)\end{bmatrix}
$$
where $K$ is given in\eqref{6.10}. Since $K$ is invertible, it follows by Step 1 that the polynomial
\begin{equation}
\widetilde f(z)=\begin{bmatrix}1 & z & \ldots &
z^{k_1+k_2-\nu_2-1}\end{bmatrix}K^{-1}
\begin{bmatrix}I_{k_1} & 0 & 0 \\ 0 & 0& I_{k_2-\nu_2}\end{bmatrix}
\begin{bmatrix}\Delta_{\boldsymbol\ell}(\balpha_1; \, g)\\ \Delta_{\boldsymbol\ell}(\balpha_2; \,
g)\end{bmatrix}
\label{6.14}
\end{equation}
satisfies conditions \eqref{6.7}. Due to the current choice \eqref{6.13} of $c_{i,j}$, we have
$$
\Delta_{\boldsymbol\ell}(\balpha_1; \, \widetilde f)=\Delta_{\boldsymbol\ell}(\balpha_1; \,
g)\quad \mbox{and}\quad \Delta_{\boldsymbol\ell}(\balpha_2; \,
\widetilde f)=\Delta_{\boldsymbol\ell}(\balpha_2; \,  g).
$$
Hence, the polynomial $f=\widetilde f-g$ satisfies homogeneous conditions \eqref{6.7} and hence,
by Step 3, 
\begin{equation}
\Delta_{\boldsymbol\ell}(\balpha_3; \, f)=\Delta_{\boldsymbol\ell}(\balpha_3; \,
\widetilde f)-\Delta_{\boldsymbol\ell}(\balpha_3; \, g)=0.
\label{6.15}
\end{equation}
We now have from \eqref{6.15}, \eqref{6.14},
\eqref{2.12} and \eqref{2.15},
\begin{align}
\Delta_{\boldsymbol\ell}(\balpha_3; g)=&\Delta_{\boldsymbol\ell}(\balpha_3; \widetilde{f})\notag\\
=&\begin{bmatrix}E_k & \mathcal J_{\balpha_3}E_{k_3} & \ldots & \mathcal
J_{\balpha_3}^{k_1+k_2-\nu_2-1}E_{k_3}\end{bmatrix}
K^{-1}
\begin{bmatrix}I_{k_1} & 0 & 0 \\ 0 & 0& I_{k_2-\nu_2}\end{bmatrix}
\begin{bmatrix}\Delta_{\boldsymbol\ell}(\balpha_1; \, g)\\ \Delta_{\boldsymbol\ell}(\balpha_2; \,
g)\end{bmatrix}\notag\\
=&V^{\boldsymbol\ell}_{k_1+k_2-\nu_2}(\balpha_3)K^{-1}
\begin{bmatrix}I_{k_1} & 0 & 0 \\ 0 & 0& I_{k_2-\nu_2}\end{bmatrix}
\begin{bmatrix}\Delta_{\boldsymbol\ell}(\balpha_1; \, g)\\ \Delta_{\boldsymbol\ell}(\balpha_2; \,
g)\end{bmatrix},\notag
\end{align}
which implies \eqref{6.5} with
$$
A=V^{\boldsymbol\ell}_{k_1+k_2-\nu_2}(\balpha_3)K^{-1}\begin{bmatrix}I_{k_1} \\ 0\end{bmatrix},\quad
B=V^{\boldsymbol\ell}_{k_1+k_2-\nu_2}(\balpha_3)K^{-1}\begin{bmatrix}0 & 0 \\ 0 
&I_{k_2-\nu_2}\end{bmatrix}.
$$
Since formula \eqref{6.5} holds for any polynomial $g\in\bH[z]$, we have, in particular,
\begin{equation}
\Delta_{\boldsymbol\ell}(\balpha_3; \, z^j)=A\Delta_{\boldsymbol\ell}(\balpha_1; \, z^j)
+B\Delta_{\boldsymbol\ell}(\balpha_2; \, z^j) \quad\mbox{for all}\quad j\ge 0
\label{6.16}
\end{equation}
and, since $\Delta_{\boldsymbol\ell}(\balpha_1; \, z^j)$,
$\Delta_{\boldsymbol\ell}(\balpha_2; \, z^j)$, $\Delta_{\boldsymbol\ell}(\balpha_3; \, z^j)$ are the 
$j$-th columns in the matrices
$V_m^{\boldsymbol\ell}(\balpha_1)$, $V_m^{\boldsymbol\ell}(\balpha_2)$,
$V_m^{\boldsymbol\ell}(\balpha_3)$ respectively, equality \eqref{6.6} follows.\qed

\smallskip

In terms of Section 3, the last theorem can be reformulated as follows.

\medskip

\begin{remark}
{\rm Given two spherical chains $\balpha_1,\balpha_2$ in the same
conjugacy class  $S\subset\bH$, let $P_{\balpha_1}$ and $P_{\balpha_2}$ be the 
associated indecomposable polynomials \eqref{3.6}. Then for any spherical chain
$\balpha_3$ whose associated polynomial $P_{\balpha_3}$ is a left divisor of the 
${\bf lrcm}(P_{\balpha_1}, P_{\balpha_2})$, equality \eqref{6.5} holds for 
any $f\in\bH[z]$ and matrices $A,B$ independent of $f$.}
\label{R:6.10}
\end{remark}

\smallskip

Theorem \ref{T:2.4} is a particular case of Theorem \ref{T:6.2} for which, however, 
we will present more explicit formulas for the coefficient matrices $A$ and $B$.
With a spherical chain $\balpha$, we associate the square matrices
\begin{equation}   
V_{\balpha}:=V^{\boldsymbol\ell}_{k}({\balpha})\quad \mbox{and}\quad
T_{\balpha}=V_{\balpha}^{-1}{\mathcal J}_{\balpha}^k
V_{\balpha}, \qquad \balpha=(\alpha_1,\ldots,\alpha_k),
\label{5.1}
\end{equation}
where ${\mathcal J}_{\balpha}$ is given in \eqref{2.12}. Recall that
$V_{\balpha}$ is invertible as the square upper triangular matrix with
all diagonal entries equal one.  

\medskip

\begin{theorem}
Under the assumptions of Theorem \ref{T:2.4}, equality \eqref{1.11} holds 
for every $f\in\bH[z]$ with 
\begin{equation} 
\begin{array}{ll}
A&=V_{\balpha_3}(T_{\balpha_3}-T_{\balpha_2})(T_{\balpha_1}-T_{\balpha_2})^{-1}
V_{\balpha_1}^{-1},\\ [1mm]
B&=V_{\balpha_3}(T_{\balpha_3}-T_{\balpha_1})(T_{\balpha_2}-T_{\balpha_1})^{-1}
V_{\balpha_2}^{-1},\end{array}\label{5.3}
\end{equation}
where $V_{\balpha_i}$ and $T_{\balpha_i}$ are $k\times k$ matrices
defined as in \eqref{5.1}.
\label{T:2.4a}
\end{theorem}

\smallskip
\noindent
{\bf Proof:}
 We get Theorem \ref{T:2.4} by letting $k_1=k_2=k_3=k$ and $\nu_2=0$ in 
Theorem \ref{T:6.2}. Therefore, the formula \eqref{1.11} holds with 
\begin{equation}
A=V^{\boldsymbol\ell}_{2k}(\balpha_3)K^{-1}\begin{bmatrix}I_{k} \\ 0\end{bmatrix},\quad
B=V^{\boldsymbol\ell}_{2k}(\balpha_3)K^{-1}\begin{bmatrix}0 \\ I_{k}\end{bmatrix},
\label{5.30}
\end{equation}
where, according to \eqref{6.10}, $K=V_{2k}^{\boldsymbol\ell}(\balpha_1,\balpha_2)$.
Making use of notation \eqref{2.15}, we can write 
\begin{equation}
V^{\boldsymbol\ell}_{2k}(\balpha_3)=\begin{bmatrix}V_{\balpha_3} & \mathcal 
J_{\balpha_3}^kV_{\balpha_3}\end{bmatrix}=V_{\balpha_3}\begin{bmatrix}I_k & 
T_{\balpha_3}\end{bmatrix}
\label{5.6}
\end{equation}
and similarly,
$$
V^{\boldsymbol\ell}_{2k}(\balpha_1,\balpha_2)=\begin{bmatrix}V_{\balpha_1} & 0 \\ 0 & 
V_{\balpha_2}\end{bmatrix}
\begin{bmatrix}I_k & T_{\balpha_1} \\ I_k & T_{\balpha_2}\end{bmatrix}.
$$
Since $V^{\boldsymbol\ell}_{2k}(\balpha_1,\balpha_2)$ is invertible, the right factor on 
the right side of  the latter equality is invertible and hence, $(T_{\balpha_2}-T_{\balpha_1})$ 
is invertible. We have
$$
V^{\boldsymbol\ell}_{2k}(\balpha_1,\balpha_2)^{-1}=\begin{bmatrix}T_{\balpha_2} & -T_{\balpha_1}\\
-I_k & I_k\end{bmatrix}\begin{bmatrix}(T_{\balpha_2}-T_{\balpha_1})^{-1}V_{\balpha_1}^{-1} & 0
\\ 0 & (T_{\balpha_2}-T_{\balpha_1})^{-1}V_{\balpha_2}^{-1}\end{bmatrix},
$$
which together with \eqref{5.6} leads us to
\begin{align*}
&V^{\boldsymbol\ell}_{2k}(\balpha_3)V^{\boldsymbol\ell}_{2k}(\balpha_1,\balpha_2)^{-1}\\
&=V_{\balpha_3}\begin{bmatrix}I_k & T_{\balpha_3}\end{bmatrix}
\begin{bmatrix}T_{\balpha_2} & -T_{\balpha_1}\\
-I_k & I_k\end{bmatrix}\begin{bmatrix}(T_{\balpha_2}-T_{\balpha_1})^{-1}V_{\balpha_1}^{-1} & 0
\\ 0 & (T_{\balpha_2}-T_{\balpha_1})^{-1}V_{\balpha_2}^{-1}\end{bmatrix}\\
&=V_{\balpha_3}\begin{bmatrix}(T_{\balpha_2}-T_{\balpha_3})(T_{\balpha_2}-T_{\balpha_1})^{-1}
V_{\balpha_1}^{-1}& 
(T_{\balpha_3}-T_{\balpha_1})(T_{\balpha_2}-T_{\balpha_1})^{-1}V_{\balpha_2}^{-1}\end{bmatrix}.
\end{align*}
Substituting the latter equality into \eqref{5.30}, implies \eqref{5.3}, thus completing the proof of the theorem.
\qed

\medskip

\begin{example}
{\rm Let $k=1$ and $\balpha_i=(\alpha_i)\subset S$ for $i=1,2,3$. According to \eqref{5.1} and \eqref{2.12},
$V_{\balpha_i}=1$ and $T_{\balpha_i}={\mathcal J}_{\balpha_i}=\alpha_i$, while formulas \eqref{1.10} and
\eqref{2.5} show that $\Delta_{\boldsymbol\ell}(\balpha_i;f)=[\alpha_i;f]_{\boldsymbol\ell}=f^{\bl}(\alpha_i)$. 
In this case, formula \eqref{1.11} amounts to \eqref{rep1}.}
\label{E:5.4}
\end{example}

\smallskip

It can be shown that the matrices \eqref{5.3} are lower triangular (although the matrices $V_{\balpha_i}$ is upper 
triangular and $T_{\balpha_i}$ is not triangular at all). Furthermore, the spherical chains  
$\balpha_1=(\alpha,\ldots,\alpha)$ and $\balpha_2=(\overline{\alpha},\ldots,\overline{\alpha})$ of the form \eqref{1.9}
satisfy the assumptions of Theorem \ref{T:2.4} and, on account of \eqref{2.8}, we arrive at the following conclusion.

\medskip

\begin{remark}
{\rm The left divided 
difference $[\gamma_1,\ldots,\gamma_k; \, f]_{\boldsymbol\ell}$ based on the spherical chain 
$\boldsymbol\gamma=(\gamma_1,\ldots,\gamma_k)\subset S$ is a left linear combination of $2k$ elements 
$(f^{(j)})^{\bl}(\alpha)$ and $(f^{(j)})^{\bl}(\overline{\alpha})$ ($j=0,\ldots,k-1$) with the coefficients
depending on $\boldsymbol\gamma$ and $\alpha$, where $\alpha$ is an arbitrary fixed element in $S$}.
\label{R:5.4}
\end{remark}

\medskip

\begin{remark}
{\rm If $\gamma_1=\ldots=\gamma_\nu=\alpha\neq \gamma_{\nu+1}$, then $[\gamma_1,\ldots\gamma_j; \, f]_{\boldsymbol\ell}$
equals $\frac{(f^{(j-1)})^{\bl}(\alpha)}{(j-1)!}$ for $j=1,\ldots,\nu$ (by \eqref{2.8}) and it is 
a left linear combination of $(f^{(j)})^{\bl}(\alpha)$ ($j=0,\ldots,j-1$)
and $(f^{(i)})^{\bl}(\overline{\alpha})$ ($i=0,\ldots,j-\nu-1$) for $j>\nu$. }
\label{R:5.4a}
\end{remark}

\medskip

\begin{example}
{\rm Let $\balpha_1=(\alpha_1,\alpha_1)$, $\balpha_2=(\alpha_2,\alpha_2)$ and $\balpha_3=(\alpha_1,\alpha_2)$
for two elements $\alpha_1\sim\alpha_2\in\bH$. Then 
\begin{align*}
V_{\balpha_1}=V_{\balpha_3}&=\begin{bmatrix} 1 &\alpha_1 \\ 0 & 1\end{bmatrix},\quad
V_{\balpha_2}=\begin{bmatrix} 1 &\alpha_2 \\ 0 & 1\end{bmatrix},\\
T_{\balpha_3}=V_{\balpha_3}^{-1}\mathcal J_{\balpha_3}^2V_{\balpha_3}&=
\begin{bmatrix} 1 &-\alpha_1 \\ 0 & 1\end{bmatrix}\begin{bmatrix} \alpha_1^2 & 0 \\ \alpha_2+\alpha_1 & 
\alpha_2^2\end{bmatrix}\begin{bmatrix} 1 &\alpha_1 \\ 0 & 1\end{bmatrix}
=\begin{bmatrix} -\alpha_1\alpha_2 & -\alpha_1\alpha_2^2-\alpha_1\alpha_2\alpha_1 \\
\alpha_2+\alpha_1 & \alpha_2^2+\alpha_2\alpha_1+\alpha_1^2\end{bmatrix},
\end{align*}
and similarly,
$\; T_{\balpha_1}=\begin{bmatrix} -\alpha_1^2 & -2\alpha_1^3 \\
2\alpha_1 & 3\alpha_1^2\end{bmatrix},\quad T_{\balpha_2}=\begin{bmatrix} -\alpha_2^2 & -2\alpha_2^3 \\
2\alpha_2 & 3\alpha_2^2\end{bmatrix}$.

\smallskip

By \eqref{2.8}, the equality \eqref{1.11} now takes the form 
\begin{equation}
\begin{bmatrix}f^{\bl}(\alpha_1) \\ [\alpha_1,\alpha_2; \, f]_{\boldsymbol\ell}\end{bmatrix}=
A\begin{bmatrix}f^{\bl}(\alpha_1)\\ (f^\prime)^{\bl}(\alpha_1)\end{bmatrix}+B
\begin{bmatrix}f^{\bl}(\alpha_2)\\ (f^\prime)^{\bl}(\alpha_2)\end{bmatrix}
\label{5.19}
\end{equation}
where (as quite tedious calculations of the right hand side expressions in \eqref{5.3} show)
$$
A=\begin{bmatrix}1 & 0 \\ (\overline{\alpha}_2-\alpha_2)^{-1} & 
(\alpha_2-\overline{\alpha}_2)^{-1}(\alpha_2-\overline{\alpha}_1)^{-1}
\end{bmatrix},\quad 
B=\begin{bmatrix}0 & 0 \\ (\alpha_2-\overline{\alpha}_2)^{-1} & 0\end{bmatrix}.
$$
It is readily seen that the equality of the top entries in \eqref{5.19} is trivial, while the 
comparison the bottom entries proves formula \eqref{2.10}.}
\label{E:5.3}
\end{example}

\medskip
\noindent
{\bf Proof of Theorem \ref{T:1.2}:} Given spherical chains \eqref{1.13},
let $\kappa$ be the integer defined in Theorem \ref{T:1.2}. By Corollary \ref{C:4.3}, 
the $\kappa$ leftmost columns in the matrix $V^{\boldsymbol\ell}_m(\balpha_1,\ldots,\balpha_n)$ (for 
$m\ge \kappa$) are right linearly independent. Thus, ${\rm rank} 
V^{\boldsymbol\ell}_m(\balpha_1,\ldots,\balpha_n)=m$ if $m\le \kappa$ and 
${\rm rank} V^{\boldsymbol\ell}_m(\balpha_1,\ldots,\balpha_n)\ge \kappa$ if $m>\kappa$. It remains to 
show that
\begin{equation}
{\rm rank} V^{\boldsymbol\ell}_m(\balpha_1,\ldots,\balpha_n)\le \kappa\quad\mbox{if}\quad m>\kappa.
\label{8.1}
\end{equation}
Let $S\subset \mathbb H$ be a conjugacy class that contains at least three spherical chains 
from \eqref{1.13}. Let $\balpha_i$ be a chain with the maximal length and for any other chain $\balpha_j\subset 
S$, denote by $\nu_j$ the number of leftmost elements shared by 
$\balpha_j$ and $\balpha_i$. Let $\balpha_r=(\alpha_{r,1},\ldots,\alpha_{r,k_r})$ be the chain in $S$ for which the
integer $k_j-\nu_j$ is maximally possible. In terms of associated polynomials \eqref{3.6}, we choose two 
polynomials $P_{\balpha_i}$ and $P_{\balpha_r}$ with maximally possible degree of their {\bf lrcm}. 

\smallskip

Any spherical chain $\balpha_j\subset S$ different from $\balpha_i$ and $\balpha_r$ satisfies conditions 
$$
k_j\le k_i\quad\mbox{and}\quad k_j-\nu_j \le k_r-\nu_r.
$$
By Theorem \ref{T:6.2}, equality \eqref{6.6} holds for all $m\ge 1$.
So, removing the block $V_m^{\boldsymbol\ell}(\balpha_j)$ from
$V^{\boldsymbol\ell}_m(\balpha_1,\ldots,\balpha_n)$ does not change the rank of the latter matrix.
Repeating this argument, we remove from $V^{\boldsymbol\ell}_m(\balpha_1,\ldots,\balpha_n)$ all blocks
$V_m^{\boldsymbol\ell}(\balpha_j)$ based on the chains in $S$ different from $\balpha_i$ and $\balpha_r$.
We then observe that the $\nu_r$ top rows in $V_m^{\boldsymbol\ell}(\balpha_i)$ are identical
to the corresponding rows in $V_m^{\boldsymbol\ell}(\balpha_r)$. Hence, removing the $\nu_r$ top rows 
from the block $V_m^{\boldsymbol\ell}(\balpha_r)$ does not change the rank of the matrix.
After this removal, the remaining matrix is of the same rank as the original one but contains only two 
blocks  ($V_m^{\boldsymbol\ell}(\balpha_i)$ and $\begin{bmatrix}0 & 
I_{k_r-\nu_r}\end{bmatrix}V_m^{\boldsymbol\ell}(\balpha_r)$)
associated with the elements in the class $S$. The total numbers of rows in these blocks 
equal (by Proposition \ref{P:3.4})
$$
k_i+k_r-\nu_r=\mu(S)=\deg({\bf lrcm}(P_{\balpha_j}: \; \balpha_j\subset S)).
$$
Repeating the latter procedure for each conjugacy class $S_j$ containing more than two spherical chains 
from \eqref{1.13}, we come up with the matrix of the the same rank as the original matrix 
$V^{\boldsymbol\ell}_m(\balpha_1,\ldots,\balpha_n)$ but having only $\kappa=\sum_j \mu(S_j)$ rows.
Its rank cannot exceed $\kappa$ which implies \eqref{7.1} and completes the proof of Theorem 
\ref{T:1.2}.\qed

\section{Formal power series over quaternions}
\setcounter{equation}{0}
In this section, we discuss several analogs of the preceding 
results in the context of the space $\bH[[z]]$ of formal power series over $\bH$. 
Given an $f\in\bH[[z]]$, we
denote by $f^\sharp$ is {\em conjugate power series} defined by
\begin{equation}
f^\sharp(z)=\sum_{j=0}^\infty z^j \overline{f}_j\quad\mbox{if}\quad f(z)=\sum_{j=0}^\infty z^j f_j.
\label{9.0}   
\end{equation}
The anti-linear involution $f\mapsto f^\sharp$ can be viewed as an extension of the 
quaternionic conjugation $\alpha\mapsto \overline{\alpha}$ from $\bH$ to $\bH[[z]]$.

\medskip
\noindent
{\bf 6.1. Linear independence of certain power series.} We consider the power series
\begin{equation}
{\bf k}_{\alpha}(z)=\sum_{k=0}^\infty \alpha^kz^k\quad\mbox{and}\quad
{\bf k}_{\alpha}^\sharp(z)={\bf k}_{\overline\alpha}(z)=\sum_{k=0}^\infty \overline\alpha^kz^k
\qquad (\alpha\in\bH).
\label{9.1}
\end{equation}
Given a spherical chain $\balpha=(\alpha_1,\ldots,\alpha_k)$, we define the 
infinite matrix $V_\infty^{\boldsymbol\ell}(\balpha)$ by letting $m=\infty$
in \eqref{1.2} and we use the entries from the same row in $V_\infty^{\boldsymbol\ell}(\balpha)$
to define power series
\begin{align}
f_1(z)=&\sum_{j=0}^\infty \left[\alpha_1;z^j\right]_{\boldsymbol\ell}z^j=\sum_{j=0}^\infty \alpha_1^jz^j=
{\bf k}_{\alpha_1}(z),\label{9.2}\\
f_2(z)=&\sum_{j=0}^\infty \left[\alpha_1,\alpha_2;z^j\right]_{\boldsymbol\ell}z^j=
\sum_{j=1}^\infty \bigg(\sum_{i=0}^{j-1} \alpha_2^i\alpha_1^{j-i}\bigg)z^j=
z{\bf k}_{\alpha_2}(z){\bf k}_{\alpha_1}(z),\notag\\
& \cdots \qquad\qquad  \cdots \qquad\qquad  \cdots \qquad\qquad  \cdots  \notag\\
f_k(z)=&\sum_{j=0}^\infty \left[\alpha_1,\ldots,\alpha_k;z^j\right]_{\boldsymbol\ell}z^j=
z^{k-1}{\bf k}_{\alpha_k}(z)\cdots {\bf k}_{\alpha_2}(z){\bf k}_{\alpha_1}(z).\notag
\end{align}
By \eqref{9.0}, \eqref{9.1} and \eqref{7.0}, the conjugate power series $f_j^\sharp$ are given by
\begin{align}
f_j^\sharp(z)&=z^{j-1}{\bf k}^\sharp_{\alpha_1}(z)\cdots
{\bf k}^\sharp_{\alpha_j}(z)=z^{j-1}{\bf k}_{\overline\alpha_1}(z)\cdots
{\bf k}_{\overline\alpha_j}(z)\label{9.3a}\\
&=\sum_{i=0}^\infty \overline{\left[\alpha_1,\ldots,\alpha_j;z^i\right]_{\boldsymbol\ell}}z^i
=\sum_{i=0}^\infty \left[\overline\alpha_1,\ldots,\overline\alpha_j;z^i\right]_{\bf r}z^i\quad
(j=1,\ldots,k).\notag
\end{align}
\begin{proposition}
Given spherical chains $\balpha_i=(\alpha_{i,1},\ldots,\alpha_{i,k_i})$ ($i=1,\ldots,n$), let
$\kappa$ be the integer defined in Theorem \ref{T:1.2}. Then the dimension of the left linear span
of power series
\begin{equation}
f_{i,j}(z)=z^{j-1}{\bf k}_{\alpha_{i,j}}(z)\cdots {\bf k}_{\alpha_{i,1}}(z)\quad
(i=1,\ldots,n; \, j=1,\ldots,k_i)
\label{9.3}
\end{equation}
and the dimension of the right linear span of the conjugate power series
\begin{equation}
f_{i,j}^{\sharp}(z)=z^{j-1}{\bf k}_{{\overline\alpha}_{i,1}}(z)\cdots {\bf k}_{{\overline\alpha}_{i,j}}(z)
\quad (i=1,\ldots,n; \, j=1,\ldots,k_i)
\label{9.3b}  
\end{equation}
equal $\kappa$. In particular, the series \eqref{9.3} are left linearly independent 
(the series \eqref{9.3b} are right linearly independent) if and only if the 
leading  elements $\alpha_{1,1},\ldots,\alpha_{n,1}$ are all distinct
and none three of them belong to the same conjugacy class.
\label{P:9.1}
\end{proposition}

\smallskip
\noindent

{\bf Proof:}
By \eqref{9.2}, the coefficients of the series \eqref{9.3} are the entries from the same row in 
the infinite confluent Vandermonde matrix  $V_\infty^{\boldsymbol\ell}(\balpha_1,\ldots,\balpha_n)$,
and the statements concerning the series \eqref{9.3} follows from Theorem \ref{T:1.2}. 
The statements concerning the conjugate power series are now immediate.\qed

\bigskip
\noindent
{\bf 6.2. Quaternion formal power series and their evaluations.}
Most of the preceding results invoked left and right evaluations of quaternion polynomials. 
To extend them to the setting of $\bH[[z]]$, we restrict our attention to a class of power series for 
which left and right evaluation functionals make sense.
We denote by $\mathbb B_R=\left\{\alpha\in\bH: \,  |\alpha|<R\right\}$ the open ball in $\bH$
of radius $R$ centered at the origin,  and we introduce the space 
$$
\mathcal H_R=\bigg\{f(z)=\sum_{j=0}^\infty f_jz^j: \; \limsup \sqrt[k]{|f_k|}\le \frac{1}{R}\bigg\}.
$$
\begin{remark}
The power series ${\bf k}_\alpha$ \eqref{9.1} belongs $\mathcal H_R$ with $R=|\alpha|^{-1}$.
More generally, the power series $z^{j-1}{\bf k}_{\alpha_j}(z)\cdots {\bf k}_{\alpha_1}(z)$
belongs to $\mathcal H_R$ with $R=\min\{|\alpha_1|^{-1},\ldots,|\alpha_j|^{-1}\}$.
\label{R:9.1a}
\end{remark}

\smallskip

We also note that for any $f(z)=\sum f_kz^k$ in $\mathcal H_R$ and any $\alpha\in\mathbb B_R$, the 
series $\sum_{k=0}^\infty\alpha^k f_k$ and $\sum_{k=0}^\infty f_k\alpha^k$
converge absolutely, so the evaluation formulas \eqref{1.2} (with $m=\infty$) make sense. Furthermore, 
the power series $L_\alpha f$ and $R_\alpha f$ (defined as in \eqref{1.3} but with $m=\infty$) are also in $\mathcal 
H_R$. Therefore, left and right divided differences for $f\in \mathcal H_R$ can be defined via formulas 
\eqref{2.5}, \eqref{2.7} for any elements $\alpha_1,\ldots,\alpha_k\in\mathbb B_R$.

\medskip

\begin{remark}
Representation formulas for quaternion polynomials in Theorems \ref{T:2.4}, \ref{T:6.2} and \ref{T:2.4a}
hold true for all elements $g\in\mathcal H_R$ and spherical chains 
$\balpha_1,\balpha_2,\balpha_3\in\mathbb B_R$.
\label{R:9.1}
\end{remark}

\smallskip
\noindent

{\bf Proof:} We verify (the most general) Theorem \ref{T:6.2}. For $g(z)=\sum g_jz^j$ we have from 
\eqref{6.16} by linearity,
\begin{align*}
\Delta_{\boldsymbol\ell}(\balpha_3; \, g)=\sum_{j=0}^\infty\Delta_{\boldsymbol\ell}(\balpha_3; \, z^j)g_j
&=\sum_{j=0}^\infty(A\Delta_{\boldsymbol\ell}(\balpha_1; \, z^j)
+B\Delta_{\boldsymbol\ell}(\balpha_2; \, z^j))g_j\\
&=A\Delta_{\boldsymbol\ell}(\balpha_1; \, g)
+B\Delta_{\boldsymbol\ell}(\balpha_2; \, g),
\end{align*}
where convergence of all series follows since $\balpha_1,\balpha_2,\balpha_3\in\mathbb B_R$.
\qed

\bigskip
\noindent
{\bf 6.3. Square summable formal power series.} We now consider the space 
$$
{\rm H}^2=\bigg\{h(z)=\sum_{j=0}^\infty z^jh_j: \; \|h\|_{{\rm H}^2}^2:=\sum_{j=0}^\infty
|h_j|^2<\infty\bigg\}
$$
of elements in $\bH[[z]]$ with square summable coefficients,
endowed with the left and right inner products
\begin{equation}
\langle h, \, g\rangle_{\boldsymbol\ell}=\sum_{j=0}^\infty \overline{g}_jh_j,\quad
\langle h, \, g\rangle_{\bf r}=\sum_{j=0}^\infty h_j \overline{g}_j.
\label{9.4}
\end{equation}
It is clear that  $\mathcal H_R\subset {\rm H}^2\subset \mathcal H_1$ (for each $R>1$) 
as sets. The series \eqref{9.1} belongs to ${\rm H}^2$ if and only if 
$\alpha\in\mathbb B_1$ (i.e., $|\alpha|<1$), and is of particular interest due to the following 
reproducing property:
$$
\langle h, \, k^\sharp_{\alpha}\rangle_{\boldsymbol\ell}=\sum_{j=0}^\infty \alpha^k h_k=h^{\bl}(\alpha)
\quad\mbox{for all}\quad h\in {\rm H}^2.
$$
More generally, if $\alpha_1,\ldots,\alpha_j$ belong to $\mathbb B_1$, then the power series
$$
f_j(z)=z^{j-1}{\bf k}_{\alpha_j}(z)\cdots {\bf k}_{\alpha_1}(z)
$$ 
belongs ${\rm H}^2$ (by Remark \ref{R:9.1}) and  reproduces the $j$-th divided difference
\begin{equation}
\langle h, \, f^\sharp_j\rangle_{\boldsymbol\ell}=\left[\alpha_1,\ldots,\alpha_j; 
h\right]_{\boldsymbol\ell}
\quad\mbox{for all}\quad h\in {\rm H}^2,
\label{9.5}
\end{equation}
which is verified by a straightforward power series computation.

\bigskip
\noindent
{\bf 6.4. Cauchy matrices.}
Given spherical chains $\balpha_i=(\alpha_{i,1},\ldots,\alpha_{i,k_i})\subset \mathbb B_1$
($i=1,\ldots,n$), the associated the power series \eqref{9.3} and their conjugates \eqref{9.3b} belong 
to ${\rm H}^2$. Let us denote by $P_{\mathcal F}^{\boldsymbol\ell}$ the left Gram matrix of the set 
\begin{equation}
\label{9.6}
\mathcal F=\{f_{i,j}^{\sharp}(z): \, i=1,\ldots,n; \, j=1,\ldots,k_i\}
\end{equation}
where $f_{i,j}^{\sharp}(z)=z^{j-1}{\bf k}_{{\overline\alpha}_{i,1}}(z)\cdots {\bf 
k}_{{\overline\alpha}_{i,j}}(z)$. Thus,
\begin{equation}
P_{\mathcal F}^{\boldsymbol\ell}=\left[P_{\balpha_i,\balpha_{i^\prime}}\right]_{i,i^\prime=1}^n,\quad\mbox{where}\quad
P_{\balpha_i,\balpha_{i^\prime}}=\left[\left\langle 
f_{i^\prime,j^\prime}^{\sharp}, \, f_{i,j}^{\sharp}\right\rangle_{\boldsymbol\ell}
\right]_{j=1,\ldots,k_i}^{j^\prime=1,\ldots,j^\prime}.
\label{9.7}
\end{equation}
\begin{remark}
{\rm In case $k_i=1$ for $i=1,\ldots,n$, we have $\mathcal F=\{{\bf k}_{\overline\alpha_1},\ldots,
{\bf k}_{\overline\alpha_n}\}$,
and consequently,
\begin{equation}
P_{\mathcal F}^{\boldsymbol\ell}=\left[\langle {\bf k}_{\overline\alpha_j}, \, 
{\bf k}_{\overline\alpha_i}\rangle_{\boldsymbol\ell} 
\right]_{i,j=1}^n=\bigg[\sum_{k=0}^\infty\alpha_i^k\overline{\alpha}_j^k\bigg]_{i,j=1}^n.
\label{9.8}
\end{equation}
If $\alpha_i\alpha_j=\alpha_j\alpha_i$, then $\sum_{k=0}^\infty\alpha_i^k\overline{\alpha}_j^k=
(1-\alpha_i\overline{\alpha}_j)^{-1}$. For this reason, we will refer to the matrix \eqref{9.8} as to
{\em left Cauchy matrix} and we will call the matrix \eqref{9.7} a} generalized (or confluent)
left Cauchy matrix.
\label{R:8.3}
\end{remark}
Connections between quaternion Cauchy and Vandermonde matrices were observed in \cite{abcs}.
The next theorem summarizes the confluent case.

\medskip

\begin{theorem}
Let $P_{\mathcal F}^{\boldsymbol\ell}$ be the confluent Cauchy matrix based on $n$ spherical chains 
$\balpha_i=(\alpha_{i,1},\ldots,\alpha_{i,k_i})\subset \mathbb B_1$ ($i=1,\ldots,n$) and defined as in \eqref{9.7}.
Then $P_{\mathcal F}^{\boldsymbol\ell}$ can be factored as 
\begin{equation}
P_{\mathcal 
F}^{\boldsymbol\ell}=V_\infty^{\boldsymbol\ell}(\balpha_1,\ldots,\balpha_n)V_\infty^{\boldsymbol\ell}(\balpha_1,\ldots,\balpha_n)^*,
\label{9.9}   
\end{equation}
where $V_\infty^{\boldsymbol\ell}(\balpha_1,\ldots,\balpha_n)$ is the infinite confluent Vandermonde matrix based on the same chains.
Therefore, 
\begin{enumerate}
\item ${\rm rank}P_{\mathcal F}^{\boldsymbol\ell}={\rm rank} V_\infty^{\boldsymbol\ell}(\balpha_1,\ldots,\balpha_n)$.
\item $P_{\mathcal F}^{\boldsymbol\ell}$ is positive definite if and only if the leftmost elements $\alpha_{1,1},\ldots,\alpha_{n,1}$ are all 
distinct and none three of them belong to the same conjugacy class.
\item The block $P_{i,i^\prime}$ of $P_{\mathcal F}^{\boldsymbol\ell}$ satisfies the Stein equality
\begin{equation}
P_{\balpha_i,\balpha_{i^\prime}}-\mathcal I_{\balpha_i}P_{\balpha_i,\balpha_{i^\prime}}\mathcal I_{\balpha_{i^\prime}}^* 
=E_{k_i}E_{k_{i^\prime}}^*,
\label{9.10}   
\end{equation}
where $\mathcal I_{\balpha_i}$, $\mathcal I_{\balpha_{i^\prime}}$, $E_{k_i}$, $E_{k_{i^\prime}}$ are defined via formulas 
\eqref{2.12}.
\end{enumerate}
\label{T:9.2}
\end{theorem}

\smallskip
\noindent
{\bf Proof:} Due to the block structure \eqref{9.7} and \eqref{1.14} of $P_{\mathcal F}^{\boldsymbol\ell}$ and 
$V_\infty^{\boldsymbol\ell}(\balpha_1,\ldots,\balpha_n)$, in order to prove \eqref{9.9} it suffices to verify
equalities
\begin{equation}
P_{\balpha_i,\balpha_{i^\prime}}=V_\infty^{\boldsymbol\ell}(\balpha_i)
V_\infty^{\boldsymbol\ell}(\balpha_{i^\prime})^*\quad\mbox{for}\quad i,j=1,\ldots,n.
\label{9.11}
\end{equation}
To this end we first observe that on account of formula \eqref{1.12}, 
\begin{equation}
\left[V_\infty^{\boldsymbol\ell}(\balpha_i)
V_\infty^{\boldsymbol\ell}(\balpha_{i^\prime})^*\right]_{j,j^\prime}
=\sum_{s=0}^\infty \left[\alpha_{i,1},\ldots,\alpha_{i,j};z^s\right]_{\boldsymbol\ell}\cdot
\overline{\left[\alpha_{{i^\prime},1},\ldots,\alpha_{{i^\prime},
{j^\prime}};z^s\right]_{\boldsymbol\ell}}.
\label{9.12}
\end{equation}
On the other hand, by virtue of formula \eqref{6.4},
$$
f_{i,j}^{\sharp}(z)=\sum_{s=0}^\infty 
\overline{\left[\alpha_{i,1},\ldots,\alpha_{i,j};z^s\right]_{\boldsymbol\ell}} z^s,
\quad  f_{{i^\prime},{j^\prime}}^{\sharp}(z)=\sum_{s=0}^\infty 
\overline{\left[\alpha_{{i^\prime},1},\ldots,\alpha_{{i^\prime},
{j^\prime}};z^s\right]_{\boldsymbol\ell}}z^s,
$$
which imply, by the definition \eqref{9.4} of the left inner product in ${\rm H}^2$,
$$
\left[P_{\balpha_i,\balpha_{i^\prime}}\right]_{j,j^\prime}=\left\langle
f_{{i^\prime},{j^\prime}}^{\sharp}, \, f_{i,j}^{\sharp}\right\rangle_{\boldsymbol\ell}=\sum_{s=0}^\infty
\left[\alpha_{i,1},\ldots,\alpha_{i,j};z^s\right]_{\boldsymbol\ell}\cdot
\overline{\left[\alpha_{{i^\prime},1},\ldots,\alpha_{{i^\prime},{j^\prime}};z^s\right]_{\boldsymbol\ell}}.
$$
Comparing the latter equality with \eqref{9.12} completes the verification of \eqref{9.11}.

\smallskip

Statements (1) and (2) follow from factorization \eqref{9.9} and Theorem \ref{T:1.2}. To prove \eqref{9.10},
we start with the Stein identity \eqref{2.20} (with $m=\infty$) applied to the chains $\balpha_i$ and 
$\balpha_{i^\prime}$ :
$$
V^{\boldsymbol\ell}_{\infty}(\boldsymbol\alpha_i)=\mathcal 
J_{{\boldsymbol\alpha_i}}V^{\boldsymbol\ell}_{\infty}(\boldsymbol\alpha_i)F_\infty^{\top}+E_{k_i}E_\infty^{\top},
\quad
V^{\boldsymbol\ell}_{\infty}(\boldsymbol\alpha_{i^\prime})=\mathcal
J_{{\boldsymbol\alpha_{i^\prime}}}V^{\boldsymbol\ell}_{\infty}(\boldsymbol\alpha_{i^\prime})
F_\infty^{\top}+E_{k_i^\prime}E_\infty^{\top}.
$$
Multiplying the first equality by the adjoint of the second equality on the right gives 
$$
V^{\boldsymbol\ell}_{\infty}(\boldsymbol\alpha_i)V^{\boldsymbol\ell}_{\infty}(\boldsymbol\alpha_{i^\prime})^*
=\left(\mathcal
J_{{\boldsymbol\alpha_i}}V^{\boldsymbol\ell}_{\infty}(\boldsymbol\alpha_i)F_\infty^{\top}+
E_{k_i}E_\infty^{\top}\right)\left(F_\infty V^{\boldsymbol\ell}_{\infty}(\boldsymbol\alpha_{i^\prime})^*\mathcal
J_{{\boldsymbol\alpha_{i^\prime}}}^*+E_\infty E_{k_i^\prime}^\top\right),
$$
which, on account of equalities 
$$
F_\infty^{\top}F_\infty^{\top}=I,\quad F_\infty^{\top}E_\infty=0,\quad E_\infty^{\top}E_\infty=1,
$$
simplifies to 
$$
V^{\boldsymbol\ell}_{\infty}(\boldsymbol\alpha_i)V^{\boldsymbol\ell}_{\infty}(\boldsymbol\alpha_{i^\prime})^*=
J_{{\boldsymbol\alpha_i}}V^{\boldsymbol\ell}_{\infty}(\boldsymbol\alpha_i)
V^{\boldsymbol\ell}_{\infty}(\boldsymbol\alpha_{i^\prime})^*\mathcal
J_{{\boldsymbol\alpha_{i^\prime}}}^*+E_{k_i}E_{k_i^\prime}^\top.
$$
The latter equality is equivalent to \eqref{9.10}, due to factorization \eqref{9.11}.\qed.

\smallskip

The right generalized Cauchy matrix $P_{\mathcal
F}^{\bf r}$ is defined by the formula \eqref{9.7} but with the right inner product rather the left.
The equality $P_{\mathcal F}^{\bf r}=V_\infty^{\bf r}(\balpha_1,\ldots,\balpha_n)^*V_\infty^{\bf 
r}(\balpha_1,\ldots,\balpha_n)$ is verified along the lines of the proof of Theorem \ref{T:9.2}.

\bigskip
\noindent
{\bf 6.5. Lagrange-Hermite interpolation.} To impose the Lagrange-Hermite interpolation 
conditions \eqref{4.6} we need well defined evaluation functionals. If the given spherical chains 
\eqref{1.13} are all in $\mathbb B_R$, we may formulate the interpolation problem \eqref{4.6} for power 
series from the space $\mathcal H_R$. Since the formulas \eqref{1.1} (and therefore, more general 
formulas \eqref{2.4}) hold  for any $f\in\mathcal H_R$ and $\alpha_1,\ldots,\alpha_k\in\mathbb B_R$
(see e.g., \cite{genstr,gss}), it follows that any power series $g\in\mathcal H_R$ satisfying homogeneous 
conditions \eqref{4.10} is necessarily of the form $g=G\cdot q$ for some $q\in\mathcal H_R$ where 
$G$ is the polynomial defined in \eqref{4.8}. Since the problem \eqref{4.6} is linear, we arrive at 
the following result.

\medskip

\begin{proposition}
Given spherical chains \eqref{1.13} in $\mathbb B_R$, such that their
leftmost elements $\alpha_{1,1},\ldots,\alpha_{n,1}$ are all distinct and none three of them belong to
the same conjugacy class, a power series $f$ belongs to $\mathcal H_R$ and satisfies conditions 
\eqref{4.6} if and only if
\begin{equation}
f=\widetilde f+Gh\quad \mbox{for some}\quad h\in\mathcal H_R,
\label{9.13}  
\end{equation}
where $\widetilde f$ and $G$ are the polynomials defined in \eqref{4.9}, \eqref{4.8}. 
If $R=1$, then $f$ of the form \eqref{9.13} belongs to ${\rm H}^2$ if and only if $h\in{\rm H}^2$
(since $\widetilde f$ and $G$ are both polynomials).
\label{P:9.2}
\end{proposition}

\bigskip
\noindent
{\bf 6.5. Norm-constrained interpolation in ${\rm H}^2$.} We finally consider the problem of finding
all power series $f\in{\rm H}^2$ satisfying interpolation conditions \eqref{4.6} and the additional norm
constraint $\|f\|_{{\rm H}^2}\le 1$. To relate (fairly explicitly) the norms of $f$ and the corresponding parameter $h$,
the formula \eqref{9.13} should be modified. The first step is to replace the particular solution $\widetilde f$ 
by another one which is orthogonal to $G\cdot {\rm H}^2$, the solution set of the homogeneous problem \eqref{4.10}.
To this end, observe that by the reproducing property \eqref{9.5}, interpolation conditions \eqref{4.6} can be written as
\begin{equation}
\langle f, \, f^\sharp_{i,j}\rangle_{\boldsymbol\ell}=c_{i,j} \quad\mbox{for}\quad
i=1,\ldots,n; \; j=1,\ldots,k_i,
\label{9.13a}
\end{equation}
where $f^\sharp_{i,j}$ are defined in \eqref{9.3}. Thus, $g\in{\rm H}^2$ satisfies homogeneous
conditions \eqref{4.10} (i.e., $g\in G\cdot {\rm H}^2$) if and only if it is left-orthogonal to all elements $f^\sharp_{i,j}$ 
from the set \eqref{9.6} and hence, to their right linear span
$$
\mathcal M={\rm span}_{\bf r}\left\{f_{i,j}^{\sharp}(z): \, i=1,\ldots,n; \, j=1,\ldots,k_i\right\}.
$$
In other words, $\mathcal M^\perp=G\cdot {\rm H}^2$, and thus we are looking for a power series
$f_{\rm min}\in\mathcal M$ satisfying conditions \eqref{9.13a}. Writing $f_{\rm min}$ in the form
\begin{equation}
f_{\rm min}(z)=\sum_{i=1}^n\sum_{j=1}^{k_1} f_{i,j}^\sharp(z)d_{ij}
\label{9.14}
\end{equation}
with unknown coefficients $d_{i,j}\in\bH$, we take its left inner product against all elements from \eqref{9.6}
getting (on account of interpolation conditions \eqref{9.13a}), the linear system
\begin{equation}
\sum_{i=1}^n\sum_{j=1}^{k_i} \langle f_{i,j}^\sharp, \,
f_{i^\prime,j^\prime}^\sharp\rangle_{\boldsymbol\ell}\cdot
d_{ij}=c_{i^\prime,j^\prime}\quad \mbox{for all}\quad i^\prime=1,\ldots,n; \, j^\prime=1,\ldots,k_{i^\prime}.
\label{9.15}  
\end{equation}
The matrix of this system is $P_{\mathcal F}^{\boldsymbol\ell}$ which is invertible under assumptions of
Proposition \ref{P:9.2}. Substituting the coefficients $d_{i,j}$ found from \eqref{9.15} into \eqref{9.14}
we get the formula
\begin{equation}
f_{\rm min}(z)=\begin{bmatrix}F_1(z) &\ldots & F_n(z)\end{bmatrix}(P_{\mathcal
F}^{\boldsymbol\ell})^{-1}C,\quad\mbox{where}\quad C=\sbm{C_1 \\ \vdots \\ C_n}
\label{9.16}
\end{equation}
is the column defined in \eqref{4.7} and where
$$
F_i(z)=\begin{bmatrix}f_{i,1}^\sharp(z) &\ldots & f^\sharp_{i,k_i}(z)\end{bmatrix}\quad\mbox{for}\quad
i=1,\ldots,n.
$$
It is not hard to show that for $f_{\rm min}$ of the form \eqref{9.16}, $\|f_{\rm min}\|_{{\rm H}^2}^2
=C^*(P_{\mathcal F}^{\boldsymbol\ell})^{-1}C$. 

\smallskip

The next modification concerns the second term on the right side of \eqref{9.13}, that is, the general 
${\rm H}^2$-solution of the homogeneous problem \eqref{4.10}.

\medskip 

\begin{proposition}
${\rm H}^2={\bf k}_{\beta}\cdot {\rm H}^2$ (as sets) for any $\beta\in\mathbb B_1$ and ${\bf k}_\beta$ defined as in 
\eqref{9.1}.
\label{P:6.4}
\end{proposition}

\smallskip
\noindent
{\bf Proof:} For any $h\in{\rm H}^2$, the series $g(z)=(z-\beta)h(z)$ also belongs to $h\in{\rm H}^2$ and, since 
${\bf k}_{\beta}(z)(z-\beta)=1$, we have $h={\bf k}_{\beta}g$ so that $h\in{\bf k}_{\beta}\cdot {\rm H}^2$.
Thus, ${\rm H}^2\subset {\bf k}_{\beta}\cdot {\rm H}^2$. 

\smallskip

To verify the reverse inclusion, take 
$h(z)=\sum h_kz^k\in{\rm H}^2$ and let $\widetilde{h}(z)=\sum |h_k|z^k$. Then $\widetilde{h}\in\mathbb R[[z]]$ 
actually belongs to the Hardy space $H^2(\mathbb D)$ of the unit disk of $\C$ and obviously satisfies
$\|\widetilde h\|_{H^2(\mathbb D)}=\|\widetilde h\|_{{\rm H}^2}=\|h\|_{{\rm H}^2}<\infty$. Since 
the power series ${\bf k}_{|\beta|}(z)=\sum|\beta|^kz^k$ is a bounded in the closed unit disk,
\begin{align*}
\|{\bf k}_{\beta}h\|^2_{{\rm H}^2}&=\sum_{k=0}^\infty \bigg|\sum_{j=0}^{k}\beta^{k-j}h_j\bigg|^2\le
\sum_{k=0}^\infty \bigg(\sum_{j=0}^{k}|\beta_{k-j}|\cdot |h_j|\bigg)^2\\
&=\|{\bf k}_{|\beta|}\widetilde 
h\|^2_{H^2(\mathbb 
D)}\le \big(\max_{z\in\mathbb \overline{D}}|{\bf k}_{|\beta|}|\big)^2\cdot \|\widetilde h\|_{H^2(\mathbb D)}=
(1-|\beta|)^{-2}\|h\|_{{\rm H}^2}.
\end{align*}
Therefore, ${\bf k}_{\beta}h$ belongs to ${\rm H}^2$ verifying the desired inclusion and completing the proof.
\qed

\medskip

\begin{proposition}
For any monic polynomial $G\in\bH[z]$ of degree $d>0$ and having no zeros outside $\mathbb B_1$, there exist
the elements $\beta_1,\ldots,\beta_d$ (each $\beta_j$ is equivalent to one of the zeros of $G$) such that 
the power series 
\begin{equation}
\Theta=G\cdot {\bf k}_{\beta_1}\cdot {\bf k}_{\beta_2}\cdots {\bf k}_{\beta_d}
\label{9.19}
\end{equation}
has the following property: $\|\Theta h\|_{{\rm H}^2}=\|h\|_{{\rm H}^2}$ for all $h\in {\rm H}^2$.
\label{P:6:5}
\end{proposition}

\smallskip

We refer to \cite[Theorem 7.1]{bol2} for the proof and explicit construction of
$\beta_1,\ldots,\beta_d$. 
In the commutative case, $\Theta$ is just the  Blaschke product having the same zeros 
(counted with multiplicities) as $G$. It follows from Proposition \ref{P:6.4} and formula 
\eqref{9.19} that $G\cdot {\rm H}^2=\Theta\cdot {\rm H}^2$ and hence, $G$ can be replaced by $\Theta$
in the parametrization formula \eqref{9.13}. We thus arrive at the modified parametrization formula 
\begin{equation}
f=f_{\rm min}+\Theta h\quad \mbox{for some}\quad h\in{\rm H^2}
\label{9.20}  
\end{equation}
with $f_{\rm min}$ given by \eqref{9.16} and $\Theta$ constructed as above. The advantages of the 
formula \eqref{9.20} can be seen from our last theorem.

\medskip

\begin{theorem}
Under the assumption of Proposition \ref{P:9.2}, the formula \eqref{9.20} describes all $f\in{\rm H}^2$
satisfying interpolation conditions \eqref{4.6}. Moreover, the representation \eqref{9.20} is left-orthogonal
and therefore,
\begin{equation}
\|f\|^2_{{\rm H}^2}=\|f_{\rm min}\|_{{\rm H}^2}^2+\|\Theta h\|_{{\rm H}^2}^2=C^*(P_{\mathcal F}^{\boldsymbol\ell})^{-1}C+\|h\|_{{\rm 
H}^2}^2.
\label{9.21}
\end{equation}
In particular, $f_{\rm min}$ is a (unique solution) to the problem \eqref{4.6} with the minimally possible norm.
In case $C^*(P_{\mathcal F}^{\boldsymbol\ell})^{-1}C\le 1$, all solutions to the problem with $\|f\|_{{\rm H}^2}\le 1$
are given by formula \eqref{9.20} with $h$ subject to $\|h\|_{{\rm H}^2}^2\le 1-C^*(P_{\mathcal F}^{\boldsymbol\ell})^{-1}C$.
\label{T:6.5}
\end{theorem}

\bibliographystyle{amsplain}

\begin{thebibliography}{10}

\bibitem{abcs}
D.~Alpay, V.~Bolotnikov, F.~Colombo and I.~Sabadini, {\em Self-mappings of the
quaternionic unit ball: multiplier properties, Schwarz-Pick inequality,  and
Nevanlinna--Pick interpolation problem}, Indiana Univ. Math. J. {\bf 64} (2015),  151-180.

\bibitem{acs}
D.~Alpay, F.~Colombo and I.~Sabadini, {\em Pontryagin-de Branges-Rovnyak spaces of slice hyperholomorphic functions}, 
J. Anal. Math. 121 (2013), 87--125.


\bibitem{bol}
V.~Bolotnikov, {\em Polynomial interpolation over quaternions}, J. Math. Anal. Appl. {\bf 421} (2015), 
no. 1, 567--590.


\bibitem{bol2}
V.~Bolotnikov, {\em 
Zeros and factorizations of quaternion polynomials: the algorithmic approach}, Preprint.

\bibitem{brenner}
J.~L.~Brenner, {\em Matrices of quaternions},
Pacific J. Math. {\bf 1} (1951), 329--335. 


\bibitem{genstr}
G.~ Gentili and D.~C.~ Struppa, {\em  A new theory of regular functions of 
a quaternionic variable}, Adv. Math. {\bf 216} (2007), no. 1, 279--301.

\bibitem{genstr1}
G.~Gentili and D.~C.~Struppa, {\em On the multiplicity of zeroes of polynomials with quaternionic coefficients},
Milan J. Math. 76 (2008), 15--25.

\bibitem{gss}
G.~Gentili, D.~C.~Struppa and C.~Stoppato, {\em Regular functions of a quaternionic variable}, Springer 
Monographs in Mathematics. Springer, Heidelberg, 2013.



\bibitem{gm}
B.~Gordon and T.~S.~Motzkin, {\em On the zeros of polynomials over division rings},
Trans. Amer. Math. Soc., {\bf 116} (1965) 218--226, 

\bibitem{kalman}
D.~Kalman, {\em The generalized Vandermonde matrix}, Math. Mag. {\bf 57} (1984), no. 1, 
15--21.


\bibitem{lam1}
T.~Y.~Lam, {\em A general theory of Vandermonde matrices}, Exposition. Math. {\bf 4} (1986), no. 3,
193--215.


\bibitem{lamler}
T.~Y.~Lam and A.~Leroy, {\em Vandermonde and Wronskian matrices over division rings},
J. Algebra {\bf 119} (1988), no. 2, 308--336.


\bibitem{lee}
H.~C.~Lee, {\em Eigenvalues and canonical forms of matrices with quaternion coefficients},
Proc. Roy. Irish Acad. Sect. A. {\bf 52}, (1949). 253--260.

\bibitem{niven}
I.~Niven, {\em Equations in quaternions}, Amer. Math. Monthly {\bf 48} (1941), 654--661.

\bibitem{ore}
O.~Ore, {\em Theory of non-commutative polynomials}, Ann. of Math. {\bf 34} (1933), no. 3, 480--508.



\bibitem{wieg}
N.~Wiegmann, {\em Some theorems on matrices with real quaternion elements},
Canad. J. Math. {\bf 7} (1955) 191--201.


\end{thebibliography}

\end{document}